\documentclass[onefignum,onetabnum]{siamonline190516}

\usepackage{amsfonts}
\usepackage{epsfig}
\usepackage[percent]{overpic}
\usepackage{amsmath}
\numberwithin{equation}{section}
\usepackage{dsfont}

\usepackage{amssymb}
\usepackage{graphics}
\usepackage{psfrag}
\usepackage{subcaption}
\usepackage{color}
\usepackage[all]{xy}
\usepackage{makeidx}
\usepackage{enumerate}
\usepackage[letterpaper,margin=0.9in]{geometry}
\usepackage{amsopn}

\usepackage{enumitem}
\setlist[enumerate]{leftmargin=.5in}
\setlist[itemize]{leftmargin=.5in}

\newsiamremark{remark}{Remark}

\def\R{\mathbb{R}}

\def\N{\mathbb{N}}

\newcommand{\be}{\begin{equation}}
\newcommand{\ee}{\end{equation}}
\newcommand{\benn}{\begin{equation*}}
\newcommand{\eenn}{\end{equation*}}
\newcommand{\bea}{\begin{eqnarray}}
\newcommand{\eea}{\end{eqnarray}}
\newcommand{\beann}{\begin{eqnarray*}}
\newcommand{\eeann}{\end{eqnarray*}}
\newcommand{\calO}{\mathcal{O}}
\DeclareMathOperator{\sgn}{sgn}

\def\txta{{\textnormal{a}}}

\def\txtd{{\textnormal{d}}}
\def\txte{{\textnormal{e}}}
\def\txtr{{\textnormal{r}}}

\def\txtD{{\textnormal{D}}}

\def\cO{\mathcal{O}}
\def\cX{\mathcal{X}}

\def\ra{\rightarrow}
\def\I{\infty}

\headers{
     Discretizations of transcritical fast-slow PDEs
}{M. Engel, C. Kuehn}

\title{ 
     Geometric analysis of a truncated
 Galerkin discretization of fast-slow PDEs with transcritical singularities
}

\author{Maximilian Engel\thanks{Department of Mathematics, Freie Universit{\"a}t Berlin, Arnimallee 6, 14195 Berlin, Germany and University of Amsterdam, KdV Institute for Mathetmatics, 1098 XG
Amsterdam, Netherlands
  (\email{m.r.engel@uva.nl}).}
\and Christian Kuehn\thanks{
Faculty of Mathematics, Technical University 
of Munich, 85748 Garching b.~M\"unchen, Germany 
(\email{ckuehn@ma.tum.de}).}}

\ifpdf
\hypersetup{
  pdftitle={
     Geometric analysis of a truncated
 Galerkin discretization of fast-slow PDEs with transcritical singularities
  },
  pdfauthor={M. Engel, C. Kuehn}
}
\fi

\begin{document}

\maketitle

\begin{abstract}
We consider a fast-slow partial differential equation (PDE) with reaction-diffusion dynamics in the fast variable and the slow variable driven by a differential operator on a bounded domain. Assuming a transcritical normal form for the reaction term and viewing the slow variable as a dynamic bifurcation parameter, we analyze the passage through the fast subsystem bifurcation point 
   for the spectral Galerkin approximation of the PDE. 
We characterize the invariant manifolds for the finite-dimensional Galerkin ODEs 
using geometric desingularization via a blow-up analysis. In addition to the crucial approximation procedure, we also make the domain dynamic during the blow-up analysis. Finally, we elaborate in which sense our results approximate the infinite-dimensional problem.
Within our analysis, we find that the PDEs appearing in entry and exit blow-up charts are quasi-linear free boundary value problems, while in the central/scaling chart we obtain a PDE, which is often encountered in classical reaction-diffusion problems exhibiting solutions with finite-time singularities.
\end{abstract}



\section{Introduction}
\label{sec:intro}

Our starting point is the classical nonlinear parabolic reaction-diffusion PDE 
\be
\label{eq:AC}
\partial_t u=\Delta u + u^2-v^2,\qquad u=u(x,t),~(x,t)\in\Omega\times[0,T],
\ee
for a smooth bounded domain $\Omega \subset \R^d$, where $v\in \R$ is the main bifurcation parameter and $T>0$ is a given fixed time, and the nonlinearity resembles the well-known transcritical normal form from ordinary differential equations (ODEs) in the reaction term. Indeed, suppose the spatial domain $\Omega=[-a,a]\subset \R$ is an interval for some fixed $a>0$, and the PDE satisfies homogeneous Neumann boundary conditions $\partial_x u(\pm a,t)=0$. Then there are two branches of constant solutions $\{u=v\}$ and $\{u=-v\}$. The local dynamics near the bifurcation point $(u,v)=(0,0)$ are well-studied, e.g., see the classical results by Crandall/Rabinowitz~\cite{CrandallRabinowitz,CrandallRabinowitz1}. There is a bifurcation at $v=0$ with a single real eigenvalue crossing the imaginary axis~\cite{Kielhoefer,KuehnBook1} and an exchange-of-stability takes place. This situation can be analyzed locally rigorously using several different methods such as Lyapunov-Schmidt reduction~\cite{Kielhoefer} or center manifold theory~\cite{HaragusIooss}. 
Yet, a comprehensive modelling approach needs to address a slowly varying parameter $v$ not only in the infinitely slow, i.e.~static limit, but also in its dynamic passage across the bifurcation point;
    given a finite amount of data, the observable (or measurable) parameter speed will typically be within a bounded range for any described natural or social system, considering phenomena such as tipping and critical transitions in ecology \cite{RinaldiMuratori}, atmospheric science \cite{Stommel}, neuroscience \cite{CornwellJones, FitzHugh} and many other applications.
Hence, a crucial case is the extension to the dynamic fast-slow bifurcation problem  
\be
\label{eq:mainPDE_fasttime}
\begin{array}{lcl}
\partial_t u&=&\partial_x^2 u + u^2-v^2 + \mu \varepsilon + H^u(u,v, \varepsilon),\\
\partial_t v&=& \varepsilon (Lv + 1 + H^v(u,v,\varepsilon)), 
\end{array} \quad u=u(x,t),~v=v(x,t),~(x,t)\in[-a,a]\times[0,T],
\ee 
where $L$ is a differential operator, 
     $\mu \in \mathbb R$ is a constant, and the terms
$H^u$ and $H^v$ are bounded smooth maps depending on higher polynomial orders of $u,v$ and $\varepsilon$ as $u, v, \varepsilon \to 0$, i.e.~representing higher order reaction terms whose effects can be controlled close to the origin.
We still assume homogeneous Neumann boundary conditions, and $\varepsilon>0$ is small so that $v$ is formally a slow variable. 
The change to the slow time scale $\tau = \varepsilon t$ gives
\be
\label{eq:mainPDE_slowtime}
\begin{array}{rcl}
\varepsilon \partial_{\tau} u&=&\partial_x^2 u + u^2-v^2 + \mu \varepsilon + H^u(u,v, \varepsilon),\\
\partial_{\tau} v&=&  (Lv + 1 + H^v(u,v,\varepsilon)).
\end{array}
\ee
    Models of the form~\eqref{eq:mainPDE_fasttime}-\eqref{eq:mainPDE_slowtime} are of relevance when studying problems in pattern formation, particularly in mathematical biology~\cite{GiererMeinhardt,Klausmeier1,Schakenberg}. In some applications, one also considers simpler variants with $L=0$, e.g., in mathematical neuroscience this is common practice, where the $u$-variable may represent voltage/current in a neuron and the $v$-variables represent (chemical) gating variables~\cite{ErmentroutTerman}. 
    Yet, physically also the gating variables do have spatial variations. Due to simplification for modelling reasons present in the ground-breaking work of Hodgkin-Huxley and FitzHugh-Nagumo \cite{FitzHugh, HodgkinHuxley, Nagumo}, which then propagated throughout the literature, such additional spatial variations in the $v$-variables had been omitted in most of the literature. Recently suitable extensions for $L\neq 0$ have emerged such as the doubly-diffusive FitzHugh-Nagumo equation~\cite{CornwellJones}. Another example, with a similar history, is a spatial version of oceanic circulation box models, such as the Stommel model for the Atlantic Meridional Overturning Circulation (see~\cite{HummelKuehn, Stommel}).

\subsection*{Background}

Note that in the singular limit $\varepsilon = 0$ in equation~\eqref{eq:mainPDE_fasttime}, the steady state homogeneous solutions are to leading-order near $(u,v)=(0,0)$ formally still given by $ u \equiv v, u \equiv -v$. If $L\equiv 0$, then it is well-known that classical finite-dimensional methods and results apply as already demonstrated at the end of the 1990s by Butuzov, Nefedov and Schneider~\cite{ButuzovNefedovSchneider1,ButuzovNefedovSchneider2}, 
    using the method of asymptotic upper and lower solutions. Here, the PDE solutions are controlled before and after the exchange of stability at the origin, by comparing to the reaction ODE without the diffusion operator $\partial_x^2$, i.e.~comparing to homogenous solutions, and using perturbation expansions for sufficiently small $\varepsilon$.
Also from a geometric viewpoint one can combine existing techniques for the case $L=0$. Using the fact that the slow variable $v$ only gives a bounded semigroup perturbation on the linearized level allows one to follow results by Bates~\cite{BatesLuZheng1,BatesLuZheng2} and then employ center manifold theory~\cite{HaragusIooss}; for an example of this strategy for the fold point see~\cite{Avitabileetal1}.
The intrinsic mathematical reason, why classical techniques work better for $L=0$, is that, in this case, the slow variables do not contain any hidden fast components generated by the differential operator $L$. Hence, for $L=0$, it is easier to find and analyze a reduced, purely slow dynamical system on a normally hyperbolic slow manifold and then to check, how the interaction between fast and slow directions occurs in regions where normal hyperbolicity is lost. Yet, the ultimate goal would be to lift the most powerful available geometric methods for ODEs directly to the full PDE setting.

For fast-slow ODEs with transcritical singularity, we recall the pivotal work by Lebovitz and Schaar \cite{LebovitzSchaar1975} who continued slow manifolds around the singularity for the case $\mu < 1$, i.e.~along stable branches of the critical manifold, using direct asymptotic estimates in a suitable rescaling. Transcritical planar fast-slow dynamics were later analyzed in great detail, for any $\mu \in \mathbb R$, by Krupa and Szmolyan as part of a series of works on the single eigenvalue crossing problem in fast-slow ODEs \cite{ KruSzm3,KruSzm4, KruSzm2}, using geometric desingularization via the blow-up method. This method has been pioneered in multiple time scale systems by Dumortier and Roussarie~\cite{DumortierRoussarie} in the middle of the 1990s. So far, trying to lift the blow-up method to the PDE \eqref{eq:mainPDE_fasttime}-\eqref{eq:mainPDE_slowtime} has remained elusive due to considerable conceptual and technical obstacles. Already the case $L=\partial_x^2 v$ is challenging and is likely to provide a key step to treat very general differential operators in the slow variable. In this work, we provide a direct variant of the blow-up method for this case. 

The first main idea for analyzing the situation with $L=\partial_x^2$ is to use a spectral Galerkin reduction~\cite{Evans} and study a spatial semi-discretization in a larger function space. In more detail, let $\{\lambda_k\}_{k=1}^\I$ denote the eigenvalues (ordered in increasing absolute value) and $\{e_k=e_k(x)\}_{k=1}^\I$ the eigenfunctions of the Laplacian $\partial_x^2$ with Neumann boundary conditions. We denote by $\langle\cdot,\cdot\rangle$ the inner product on $L^2([-a,a])$. Then setting 
\benn
u(x,t)=\sum_{k=1}^\I u_k(t) e_k(x),\qquad v(x,t)=\sum_{k=1}^\I v_k(t) e_k(x)
\eenn 
gives, upon taking the inner product of~\eqref{eq:mainPDE_fasttime} with each $e_k$, the infinite Galerkin 
ODE system
\be
\label{eq:mainPDE_Galerkin}
\begin{array}{lcl}
\partial_t u_k&=&\lambda_k u_k + \langle u^2-v^2,e_k\rangle + \varepsilon \langle \mu ,e_k\rangle + \langle H^u(u,v, \varepsilon),e_k\rangle ,\\
\partial_t v_k&=& \varepsilon (\lambda_k v_k +\langle 1+H^v(u,v,\varepsilon),e_k\rangle). 
\end{array} 
\ee
Infinite-dimensional ODEs play a key role in nonlinear dynamics of PDEs~\cite{Henry}, e.g., in the context of amplitude/modulation equations~\cite{SchneiderUecker}. Yet, the fast-slow case is poorly understood currently as one essentially has to deal with two infinite-dimensional ODEs simultaneously due to to the scale separation. The classical strategy in PDEs used to study Galerkin systems is to first use a truncation $u_k=0=v_k$ for all $k > k_0$ and a fixed $k_0\in \N$. Then one studies the resulting finite ODE system~\cite{Evans,SchneiderUecker} with the dynamical understanding that near many common bifurcation points, only the projection onto the lowest order modes $e_k$ with $k\leq k_*$ for some intrinsic fixed critical $k_*$ should be relevant. 

Now, the crucial insight of the related work \cite{HummelKuehn} concerns the fact that for fixed $\varepsilon$ in equation~\eqref{eq:mainPDE_Galerkin} there is a $k_0(\varepsilon)$ such that $\left|\lambda_k \right|> \frac{1}{\varepsilon}$ for all $k > k_0(\varepsilon)$; hence, these higher modes as part of the ``slow'' variable have to be counted as fast modes. 
In other words, depending on the fixed $\varepsilon$, the parametrization of a normally hyperbolic slow manifold has to account for this separation by an additional parameter $\zeta$, in our specific Galerkin case directly corresponding with $k_0(\varepsilon)$. Such manifolds $S_{ \varepsilon, \zeta}$, and their specific Galerkin approximations $S_{\txta, \varepsilon}^{-,k_0}$ (see \cite{EngelHummelKuehn}), 
will be explicitly described in the following before passing the singularity at $(u,v)=(0,0)$.

 \subsection*{Main result}   
In order to deal with the loss of hyperbolicity at $u=0=v$, we employ the blow-up method at each Galerkin level. More precisely, we can still carry out a blow-up transformation for each truncated Galerkin system where $u_k=0=v_k$ for all $k > k_0$.
This leads to an infinite sequence of blown-up problems with growing dimensions of the spherical blown-up space, which converges to a Banach manifold~\cite{Lang1} as $k_0\ra +\I$, in this case an infinite-dimensional sphere~\cite{AbrahamMarsdemRatiu} 
corresponding with an appropriate $L^2$ space. The dynamics on this limiting object of the blown-up phase spaces are very challenging to describe:
Considering the intricacies of infinite-dimensional fast-slow systems as discussed in the previous paragraph,
such a sequence $k_0 \to \infty$ will generally entail $\varepsilon \to 0$ for the corresponding slow manifolds to be tracked and, hence, not necessarily give perturbation results. We will discuss in what sense first order approximations of $S_{\txta, \varepsilon}^{-,k_0}$, and its counterparts in blow-up coordinates, satisfy a fully infinite-dimensional description for our specific example.

In the following, 
we are studying the entire family of truncated Galerkin systems~\eqref{eq:mainPDE_Galerkin}, using ODE blow-up methods,
yielding the tracking of slow manifolds $S_{\txta, \varepsilon}^{-,k_0}$ around the non-hyperbolic singularity.
An important ingredient in this process is also to view the parameter $a\in\R$ specifying the domain as a dynamic variable within the blow-up. 
    This derives from the fact that the parameter $a$ enters the Galerkin ODEs \eqref{eq:mainPDE_Galerkin} via the eigenvalues $\lambda_k$. Hence, this parameter now appears in addition to $\varepsilon$ such that finding an appropriate desingularization has to involve control of $a$.
On a related note, we also identify within each chart of the blow-up construction related PDEs of independent interest such as a free boundary problem in the entry and exit charts, while recovering a more traditional, yet nonlocal, scalar reaction-diffusion problem with polynomial nonlinearity in the scaling chart. 

The analysis yields the following main result:
\begin{theorem}
\label{thm:mainresult}
The attracting slow manifolds $S_{\txta, \varepsilon}^{-,k_0}$ near the origin for system~\eqref{eq:mainPDE_Galerkin}, truncated at $k_0 \in \mathbb{N}$, 
exhibit the following behaviour:
for any fixed $\mu \neq 1$, there exists $\varepsilon_0 > 0$ such that for all $\varepsilon \in (0, \varepsilon_0]$:
\begin{itemize}
\item (exchange-of-stability) If $\mu < 1$, the manifold $S_{\txta, \varepsilon}^{-,k_0}$ locally extends around the origin inside a small tubular neighbourhood of the attracting critical manifold $\{(u,v) : u=-v < 0\}$.
\item (jump case) If $\mu > 1$, the manifold $S_{\txta, \varepsilon}^{-,k_0}$ locally extends around the origin inside a small tubular neighbourhood of the set $\{(u,v) \,:\, u> 0, \ v = 0\}$.
\end{itemize}
\end{theorem} 
The immediate corollary of Theorem~\ref{thm:mainresult} is that solutions of the reaction-diffusion PDE locally starting near $\{u=v,v<0\}$ first approach an attracting manifold $S_{\varepsilon, \zeta}$ for $\zeta$ corresponding with $k_0$ (see \cite{EngelHummelKuehn}). 
Then two cases occur depending upon $\mu$: either there is an exchange-of-stability and the Galerkin approximation $S_{\txta, \varepsilon}^{-,k_0}$ of $S_{ \varepsilon, \zeta}$ extends towards the other locally attracting part of the critical manifold, or we obtain a jump away from the bifurcation point along a fast fiber. The special case $\mu=1$ is naturally going to correspond to canards but we shall not treat this higher-codimension case here; see~\cite{BenoitCallotDienerDiener,Benoit2,DumortierRoussarie,Wechselberger1} for more details regarding canards in ODEs. 
    Note that we do not obtain an extension around the singularity of $S_{\varepsilon, \zeta}$ itself. Such a result is beyond the scope of this work for now, as it would involve a double limit of $\varepsilon \to 0$ and $k_0 \to \infty$ whose difficulties together with non-rigorous limiting approaches we discuss throughout the mansuscript.
    However, as individual trajectories of the truncated Galerkin system are close to those of the PDE for sufficiently small initial conditions, we effectively obtain the extension of ``approximately invariant slow manifolds" in form of $S_{\txta, \varepsilon}^{-,k_0}$.

\subsection*{Related recent work}
Our result in Theorem~\ref{thm:mainresult} provides a generalization of the ODE situation~\cite{KruSzm4}, while the proof technique is 
  the first PDE-related geometric desingularization, taking into account an arbitrary number of Galerkin modes.  
As such, this new method-of-proof, and not the result per-se, forms the core contribution of our work in combination with \cite{Engeletal2022}. In fact, we expect that our technique can become a blueprint for using geometric desingularization for infinite-dimensional dynamical systems as discussed in Section~\ref{sec:conclusion}.

The results on slow manifolds in the infinite-dimensional setting \cite{HummelKuehn} and their relation to Galerkin manifolds \cite{EngelHummelKuehn} build the basis for the analysis before the the singularity is reached. Then the blow-up of the Galerkin ODEs is needed to understand the transition around the singularity, yielding new challenges compared to the pure ODE setting. The fact that we deal with a domain length $a$ is solved in our case by including this parameter into the blow-up transformation, involving an additional dynamical direction in the blow-up charts. This differs crucially from \cite{Engeletal2022} where another $\varepsilon$-dependent rescaling is used to continue treating $a$ as a parameter. Both ways yield the desired results with different side effects: the approach chosen here naturally links the blow-up Galerkin ODEs to PDEs with moving boundaries.
Note that our treatment of the transcritical situation involves a case distinction depending on $\mu$, in contrast to a fold problem.
In this paper, we have additionally provided an expansion of the Galerkin-center (or slow) manifolds, yielding a detailed insight into the Galerkin dynamics and enabling a perspective on potential limits of these expansions. 

\subsection*{Structure of the paper}
The paper is structured as follows. Section~\ref{sec:beforeblowup} analyzes system~\eqref{eq:mainPDE_fasttime} and the corresponding system~\eqref{eq:mainPDE_Galerkin} around the homogenous solution $\{u=v, v<0,\varepsilon=0\}$ near the origin but in the normally hyperbolic regime using established theory to obtain suitable slow manifolds before the fast subsystem bifurcation. Before proceeding to the blow-up regime near the bifurcation, we consider in Section~\ref{sec:smallk0} the Galerkin approximation for small $\varepsilon$ by approximating the center-stable, or slow, manifolds for each $k_0 > 0$, and consider a (non-rigorous) first-order-in-$\varepsilon$ approximation in the limit $k_0 \to \infty$ (which will actually not be needed for the main blow-up result). In Section~\ref{sec:blowup}, we discuss the main dynamical results in the three different charts of an appropriately chosen blow-up manifold when normal hyperbolicity is lost, again for each fixed $k_0 > 0$, allowing us to extend the approximating slow manifolds around the non-hyperbolic origin. Furthermore, we investigate the correspondence of the Galerkin ODEs in the coordinates of the charts with rescaled PDEs, related to~\eqref{eq:mainPDE_fasttime}. Finally, we use the results of the blow-up to prove a more detailed version of Theorem~\ref{thm:mainresult} in Section~\ref{sec:mainresult}. Then we conclude with a short summary and outlook in Section~\ref{sec:conclusion}.

\section{Analysis of the Galerkin system before blow-up} 
\label{sec:beforeblowup}

\subsection{A Remark on Solution Theory} \label{sec:solution_theory}

The local-in-time solution theory of the PDE
\be
\label{eq:mainPDE_fasttimePDEview}
\begin{array}{lcl}
\partial_t u&=& \partial_x^2 u + u^2-v^2 + \mu \varepsilon + H^u(u,v, \varepsilon),\\
\partial_t v&=& \varepsilon (\partial_x^2v +1+H^v(u,v,\varepsilon)), 
\end{array} 
\ee 
for $u=u(x,t)$, $v=v(x,t)$, $(x,t)\in[-a,a]\times[0,T]$, with initial conditions $u(x,0)=u_0(x), v(x,0)= v_0(x)$ of sufficiently small $L^2$-norm, and in an appropriate space and Neumann boundary conditions
\be
\label{eq:mainPDE_fasttimePDEview1}
\partial_x u(\pm a,t)=0,\quad \partial_x v(\pm a,t)=0,~\forall t\in[0,T],
\ee
is relatively standard, and we briefly recall the relevant results in terms of classical solutions. We are interested in the case, when $H^u$ and $H^v$ are bounded smooth maps with 
$$H^u(u,v, \varepsilon)=\mathcal{O} \left(u^3, u^2 v, v^2 u, v^3, \varepsilon u, \varepsilon v, \varepsilon^2\right)$$
 and $H^v(u,v,\varepsilon)=\cO(u^{2},uv,v^{2},\varepsilon)$ as $u,v,\varepsilon\ra 0$. 
 
We may write \eqref{eq:mainPDE_fasttimePDEview} as 
\[
	w_t = A w + F(w), \quad\text{with } w(0)= w_0,
\]
where $w = (u,v)^T$, $w_0 = (u_0, v_0)^T$, $F(w) = ( u^2 - v^2 + \mu \varepsilon + H^u(u, v, \varepsilon), \varepsilon + \varepsilon H^v(u, v, \varepsilon))^T$, and 
\[
	Aw =\begin{pmatrix} 
	\partial_x^2 u & 0 \\
	0 & \varepsilon \partial_x^2 v
	\end{pmatrix}, \quad\text{with } 
	\mathcal D(A) = \{w \in  H^2(-a,a)^2 \; : \; \partial_x u = 0 = \partial_x v \; \text{ at }x=\mp a \}.  
\]
We have that $F(w)$ is locally Lipschitz continuous on $Z^\alpha = \mathcal D(A^\alpha)$ for
$1/4<\alpha <1$. Moreover, the operator $A$ is sectorial and a generator of an analytic semigroup on
$Z=L^2(-a,a)^2$. Thus, for $w_0 \in Z^\alpha$, there exists a unique local-in-time solution $w \in
C([0, t_\ast); Z^\alpha) \cap C^1((0, t_\ast); Z)$, with $w(t) \in \mathcal D(A)$, to
\eqref{eq:mainPDE_fasttimePDEview} for some $t_\ast >0$; see e.g.~\cite{Henry}. The quadratic nonlinearity in
\eqref{eq:mainPDE_fasttimePDEview}  implies a potential finite-time blowup of solutions to
\eqref{eq:mainPDE_fasttimePDEview}; cf.~e.g.~\cite{Ball}. However, simple estimates show that, for initial
$u_0 <0$ and $v_0>0$, a solution of \eqref{eq:mainPDE_fasttimePDEview} exists for $t>0$ such that $u(t)\leq 0$ and $v(t)\geq 0$.

 

Since we are just interested in the passage near $(u,v)=(0,0)$ of near-constant initial data near the critical manifold $\{u=v,v<0\}$, we could even assume without loss of generality that the higher-order terms $H^u$ and $H^v$ are chosen such that the problem is globally dissipative, when suitable norms of $u,v$ become large. This would guarantee global-in-time existence ($T=+\I$) via standard methods~\cite{Robinson1,Temam}. Yet, these technical aspects would not influence the main problem to find a geometric description of the flow via the blow-up method, so we shall not focus on them here and work with sufficiently regular and bounded solutions near the origin $(u,v)=(0,0)$.

\subsection{Explicit form of the Galerkin ODEs} \label{sec:explicit_form_of_Galerkin_ODEs}

Next, we want to implement the blow-up method. A geometric approach is considerably more elegant in a Hilbert space setting and since all bounded smooth solutions certainly satisfy $u(t),v(t)\in L^2([-a,a])$ for all $t\in[0,T]$, we carry out the geometric analysis in $L^2=L^2([-a,a])$.
Consider equation~\eqref{eq:mainPDE_fasttime} on the interval $[-a,a]$ with Neumann boundary conditions such that the eigenfunctions and eigenvalues of the Laplacian $\partial_{x}^2$ are given as
\begin{equation} \label{eigenfct:Neumann}
e_{k+1}(x) = \sqrt{\frac{1}{a}} \cos \left( \frac{\pi k (x+a)}{2a} \right) + \delta_{k0} \frac{1 - \sqrt{2}}{\sqrt{2a}}, \ \lambda_{k+1} = -\left( \frac{\pi}{2a} k \right)^2, \ k \geq 0,
\end{equation} 
where $\delta_{k0}=1$ if $k=0$ and $\delta_{k0}=0$ otherwise. In the following, we derive equation~\eqref{eq:mainPDE_Galerkin} for truncation sizes $k_0$ such that $u_k=0=v_k$ for all $k > k_0$, where we recall that\begin{align*}
H^u(u,v, \varepsilon) &= \mathcal{O} \left(u^3, u^2v, v^2u, v^3, \varepsilon u, \varepsilon v, \varepsilon^2\right), \\
H^v(u,v, \varepsilon) &= \mathcal{O} \left(u^{2},uv, v^{2}, \varepsilon\right).
\end{align*}
We obtain the following system of ODEs for any fixed $k_0 \geq 1$:
\begin{proposition} \label{truncation_wayofwriting}
Consider equation~\eqref{eq:mainPDE_fasttime} on the interval $[-a,a]$, where $L=\partial_x^2$ and $\mu \in \mathbb{R}$, such that the eigenfunctions and eigenvalues of the Laplacian with Neumann boundary conditions $\partial_x e_k(\pm a,t)=0$ are given in equation~\eqref{eigenfct:Neumann}. 
Then the system of infinite Galerkin ODEs~\eqref{eq:mainPDE_Galerkin}, truncated at $k_0 \geq 1$, can be written as
\begin{subequations}
\label{eq:k0_general}
\begin{align}
\partial_t u_1&= u_1^2 - v_1^2 + 2a \varepsilon \mu +  \sum_{j=2}^{k_0} \left(u_j^2 - v_j^2 \right) + H_{1}^u,\\
\partial_t v_1&= \varepsilon(2a + H_1^v) , \\
\partial_t u_k&= \hat \lambda_k u_k + 2 (u_k u_1  -  v_k v_1) + \sqrt{2}\sum_{i,j=2}^{k_0} \alpha_{i,j}^k (u_i u_j - v_i v_j) + H_k^u, \ 2 \leq k\leq k_0,\\
\partial_t v_k&= \varepsilon (\hat \lambda_k v_k +  H_k^v), \ 2 \leq k \leq k_0, \\
\partial_t \varepsilon &= 0,
\end{align} 
\end{subequations}
where $\hat \lambda_k = \sqrt{2a} \lambda_k $,
$\alpha_{i,j}^k \in [0, 1]$ with $\alpha_{i,j}^k \neq 0$ if and only if $i+j-k =1 \lor k-\left|i-j\right| =1$ and for all $1\leq k \leq k_0$
\begin{align*}
H_k^u &= \mathcal{O} \left(u_i^3, u_i^2v_j, v_j^2u_i, v_j^3, \varepsilon u_i, \varepsilon v_j, \varepsilon^2\right), \ i,j \in \{1, \dots, k_0\},\\
H_k^v &= \mathcal{O} \left( u_i^{2}, u_i v_j,  v_j^{2}, \varepsilon\right), \ i,j \in \{1, \dots, k_0\}.
\end{align*}
\end{proposition}

\begin{proof}
Firstly, observe that
\begin{align*}
&\left\langle \left( \sum_{j=1}^{k_0} u_j e_j \right)^2 - \left( \sum_{j=1}^{k_0} v_j e_j \right)^2, e_k\right\rangle \\ &= \sum_{j=1}^{k_0} u_j^2 \langle e_j^2, e_k\rangle + \sum_{j,i=1, j\neq i}^{k_0} u_j u_i \langle e_j e_i, e_k \rangle - \left( \sum_{j=1}^{k_0} v_j^2 \langle e_j^2, e_k\rangle + \sum_{j,i=1, j\neq i}^{k_0} v_j v_i \langle e_j e_i, e_k \rangle \right).
\end{align*}
Secondly, we can calculate
\begin{align*}
\lambda_1 &= 0, \ \langle 1, e_1 \rangle = \sqrt{2a}, \\
\langle e_k e_l, e_1 \rangle &= \langle e_1 e_l, e_k \rangle = \frac{1}{\sqrt{2a}}  \delta_{k,l}, \quad k,l \geq 1.
\end{align*}
Hence, we obtain
\begin{align*}
\partial_t u_1&= \frac{1}{\sqrt{2a}} \left(u_1^2 - v_1^2\right) + \sqrt{2a} \varepsilon \mu + \frac{1}{\sqrt{2a}} \sum_{j=2}^{k_0 -1} \left(u_j^2 - v_j^2\right) + \frac{1}{\sqrt{2a}} H_1^u ,\\
\partial_t v_1&= \varepsilon \sqrt{2a} + \frac{1}{\sqrt{2a}} \varepsilon H_1^v.
\end{align*}
Additionally, we observe that, for $2 \leq k \leq k_0$,
\begin{equation*}
\langle e_k^2, e_k \rangle = 0, \quad b_k = \langle 1, e_k \rangle = 0,
\end{equation*}
and that, for $2\leq i,j,k\leq k_0$,
\begin{align*}
\langle e_i e_j, e_k \rangle &= \frac{1}{\sqrt{a}} \frac{1}{a} \int_{-a}^a \cos \left( \frac{\pi (i-1) (x+a)}{2a} \right) \cos \left( \frac{\pi (j-1) (x+a)}{2a} \right) \cos \left( \frac{\pi (k-1) (x+a)}{2a} \right)~\txtd x \\
&= \frac{2}{\sqrt{a}} \int_{0}^1 \cos \left(  (i-1)\pi x \right) \cos \left((j-1)\pi x \right) \cos \left(  (k-1)\pi x \right) ~\txtd x \\
&= \frac{1}{\sqrt{a}} \left(\int_{0}^1 \cos \left(  (i + j-2)\pi x \right) \cos \left(  (k-1)\pi x \right) ~\txtd x
+  \int_{0}^1  \cos \left((i-j)\pi x \right) \cos \left(  (k-1)\pi x \right) ~\txtd x \right)\\
&= \frac{1}{\sqrt{a}} \alpha_{i,j}^k,
\end{align*}
where $\alpha_{i,j}^k \in [0,1]$ is non-zero if and only if $i+j-2 = k-1$, i.e.~$i+j-k =1$, or $\left|i-j\right|=k-1$, i.e.~$k-\left|i-j\right|=1$. Hence, we get 
\begin{align*}
\partial_t u_k&= \lambda_k u_k + \frac{2}{\sqrt{2a}} (u_k u_1  -  v_k v_1) + \frac{1}{\sqrt{a}} \sum_{i,j=2}^{k_0} \alpha_{i,j}^k (u_i u_j - v_i v_j) + \frac{1}{\sqrt{2a}} H_k^u, \ k\geq 2,\\
\partial_t v_k&= \varepsilon \lambda_k v_k + \frac{1}{\sqrt{2a}}\varepsilon H_k^v, \ k\geq 2, \\
\partial_t \varepsilon &= 0.
\end{align*} 
Furthermore, we conduct the time change $\tilde{t} = t/\sqrt{2a}$, and then drop the tilde again, which finishes the proof.
\end{proof}

\subsection{Slow and Galerkin manifolds}
\label{sec:manifolds}

In analogy to a standard procedure for fast-slow ODEs of singular perturbation type, the first step in
our geometric analysis is to determine the critical manifold for \eqref{eq:mainPDE_fasttimePDEview}.
Considering the slow formulation of~\eqref{eq:mainPDE_fasttimePDEview}, obtained from the time rescaling
$\tau = \varepsilon t$, 
\begin{subequations}%
	\label{eq:initial-problem-slow}
	\begin{align}
		\varepsilon u_\tau &= \partial_x^2 u + u^2-v^2 + \mu \varepsilon + H^u(u,v, \varepsilon) && \text{ for } x\in (-a,a)\text{ and }\tau>0, \\
		v_\tau &= \partial_x^2 v +1+H^v(u,v,\varepsilon) && \text{ for } x\in (-a,a)\text{ and }\tau>0, \\
		& \partial_x u(\tau ,x) = 0 =  \partial_x v(\tau,x)  && \text{ for } x=\mp a \text{ and }\tau>0,
	\end{align}
\end{subequations}
and setting $\varepsilon = 0$ therein, we find that the critical manifold is given by the set
\begin{equation}
	\left\{ (u,v) : 0 = \partial_x^2 u + u^2-v^2 + \mu \varepsilon + H^u(u,v, \varepsilon),~ \partial_x u(\cdot,\mp a)=0= \partial_x v(\cdot,\mp a) \right\}.
\end{equation}
Restricting to spatially homogeneous solutions, we define the critical manifold $S_0$ as the set of functions
\begin{equation}%
	\label{eq:critical-manifold-pde}
	S_0 := \left\{ (u,v) \in \R^2 : 0 = u^2-v^2 + H^u(u, v, 0) \right\},
\end{equation}
abusing notation and identifying constant functions $u : [-a, a] \rightarrow \R$ with the value $u$ they take.
Due to our assumptions on the form of $H^u$, near the origin $(u, v) = (0,0)$ in $(u,v)$-space the set $S_0$ is given as 
\begin{equation}
	S_0 = \left\{ (u, v) \in \R^2: v = \pm u + \calO\left( u^3 \right) \right\}.
\end{equation}
Proceeding, again, as in a finite-dimensional setting, the second step in our analysis
concerns the persistence of the manifold $S_0$ for $\varepsilon$ positive and sufficiently small.
However, in an infinite-dimensional setting, the concept of ``fast'' and ``slow'' variables can be delicate, as for
any $\varepsilon > 0$, there exists $k>0$ such that $\varepsilon \lambda_k = O(1)$.
One way to address this complication is to split the slow variables $v$ into fast and slow parts,
which we explain in the context of the following Prop.~\ref{Prop:slowman}, where we focus on the lower left branch of $S_0$ whose continuation through the origin we will be interested in. 
We refer to~\cite{HummelKuehn}
for further discussion and details on this ``fast-slow splitting'' in infinite dimensions.
\begin{proposition}\label{Prop:slowman}
	Let $(u, v) \in S_0$ with $u =v < 0$. Consider any small $\zeta>0$ and  $u \le \omega_A < 0$, $\omega_f \in \R$, and $L_f > 0$ such that  $\omega_A + L_f
	<\omega_f < 0$. Then, there exist spaces $Y^\zeta_S \oplus Y^\zeta_F = L^2(-a, a)$ and a family
	of attracting slow manifolds around $(u, v)$ that are given as graphs
	\begin{equation}\label{slow_manig_1}
		S_{\varepsilon, \zeta} := \left\{ \left( h^{\varepsilon, \zeta}_X(v), 
			h^{\varepsilon,	\zeta}_{Y^\zeta_F}(v), v \right) : v \in Y_S^\zeta \right\}
	\end{equation}
	for $0 < \varepsilon < C \frac{\omega_f}{\omega_A} \zeta$ and some fixed  $C \in (0, 1)$,
	where $\left( h^{\varepsilon, \zeta}_X(v), h^{\varepsilon,	\zeta}_{Y^\zeta_F}(v)
	\right) :
	Y^\zeta_S \rightarrow H^2(-a,a) \times (Y^\zeta_F \cap H^2(-a,a))$. 
\end{proposition}
\begin{proof}
Analogous to \cite{Engeletal2022}.
\end{proof}
The construction of these subspaces is as follows:
Note that, for any $\varepsilon >0$, there exists $k>0$ such that $\varepsilon \lambda_k = \calO(1)$,
where  $\lambda_k$ are given in~\eqref{eigenfct:Neumann}. Thus, to define fast and slow variables, we need
to split the basic space $Y=L^2(-a,a)$ for $v$ into $Y = Y^\zeta_S \oplus Y^\zeta_F$, where
\begin{subequations}%
	\label{eq:Y_split}
	\begin{align}
		Y^\zeta_S &:= \operatorname{span} \left\{ e_k(x) : 0 \le k \le k_0 \right\}\quad\text{and} \\
		Y^\zeta_F &:= \overline{\operatorname{span} \left\{e_k(x) : k > k_0 \right\}}^{L^2},
	\end{align}
\end{subequations}
with $\{e_k(x)\}_{k \in \mathbb N}$ being the eigenfunctions \eqref{eigenfct:Neumann} corresponding to the eigenvalues $\{\lambda_k\}_{k \in \mathbb N}$ and $k_0 \in \mathbb N$ satisfying
\begin{equation}%
	\label{eq:zeta-and-k0}
	- \frac{{(k_0 + 1)}^2 \pi^2}{4a^2} < \zeta^{-1} \omega_A \le - \frac{k_0^2 \pi^2}{4 a^2}, 
\end{equation}
for  given $\zeta>0$ and $\omega_A$ given by Proposition~\ref{Prop:slowman}. Next, for given $\zeta>0$, we also split the space $X=L^2(-a,a)$ into $X= X^\zeta_S \oplus X^\zeta_F$, where $X^\zeta_S$ and $X^\zeta_F$ are defined in the same manner as $Y^\zeta_S$ and $Y^\zeta_F$, see~\eqref{eq:Y_split}. Then, the truncation of the Galerkin system in~\eqref{eq:mainPDE_Galerkin} at $k_0$, which is related to $\zeta$ via~\eqref{eq:zeta-and-k0}, gives the projection
of~\eqref{eq:mainPDE_fasttimePDEview} onto $\big( X^\zeta_S, Y^\zeta_S \big)$. Thus, we obtain a family of so-called Galerkin manifolds
\begin{equation}%
\label{eq:galerking_manif}
	G_{\varepsilon, \zeta} := \left\{ \left( h_G^{\varepsilon, \zeta} (v), v \right) : v \in
		Y^\zeta_S \right\}
\end{equation}
for a function $h^{\varepsilon, \zeta}_G : Y^\zeta_S \rightarrow X^\zeta_S$. 

\begin{proposition}\label{Prop:slow_Galerkin}
	For $0 < \varepsilon < C \frac{\omega_f}{\omega_A} \zeta$ and some fixed $C \in (0, 1)$, where $\zeta$, $\omega_A$ and $\omega_f$ are as in Proposition~\ref{Prop:slowman}, the
	following estimate holds:
	\begin{equation}%
		\left\Vert h^{\varepsilon, \zeta}_X(v) - h_G^{\varepsilon, \zeta} (v) \right\Vert_{H^{2}} + \left\Vert h^{\varepsilon,
		\zeta}_{Y^\zeta_F}(v) \right\Vert_{H^{2}} \le \tilde C \left( \frac{4a^2}{\pi^2 (2 k_0 + 1)} + \zeta
		\right) \Vert v \Vert_{H^{2}}.
	\end{equation}
	In particular, using the relation between $\zeta$ and $k_0$ in~\eqref{eq:zeta-and-k0}, we have
	\begin{equation}%
		\left\Vert h^{\varepsilon, \zeta}_X(v) - h_G^{\varepsilon, \zeta} (v) \right\Vert_{H^{2}} + \left\Vert h^{\varepsilon,
		\zeta}_{Y^\zeta_F}(v) \right\Vert_{H^{2}} \le \tilde C \frac{1}{k_0} \Vert v \Vert_{H^{2}}.
	\end{equation}
\end{proposition}
\begin{proof}
The proof follows the same steps as in~\cite{EngelHummelKuehn}.
\end{proof}

\begin{remark} Note that $k_0 \to \infty$ corresponds to $\zeta \to 0$ which, due to the relation $0<\varepsilon <
C \frac{\omega_f}{\omega_A} \zeta$, see Propositions~\ref{Prop:slowman} and~\ref{Prop:slow_Galerkin}, implies also $\varepsilon \to 0$ when $k_0
\to \infty$. Hence, the limit of the Galerkin manifolds $G_{\varepsilon, \zeta}$ as $k_0 \to \infty$
cannot, in general, be guaranteed for all $0 < \varepsilon < \varepsilon_0$ with an independent
upper bound $\varepsilon_0$. Thus, we perform the following analysis for $0 < \varepsilon < \varepsilon_0$, with $\varepsilon_0$ sufficiently
small, and $k_0$ arbitrarily large, but fixed. 
\end{remark}

The outlined situation requires two main steps: (1) in the region where we do have sufficient hyperbolicity and linear stability, we have to construct for our Galerkin system a suitable family of manifolds $S^{-,k_0}_{\txta,\varepsilon}$, which are sufficient to approximate all solutions starting near the homogeneous attracting branch of the critical manifold, and (2) we have to use a blow-up to propagate points inside $S^{-,k_0}_{\txta,\varepsilon}$ through the singularity. We shall proceed with this program in Section~\ref{sec:blowup} for the truncated Galerkin system. The reader only interested in the blow-up approach can skip ahead without loss continuity to Section~\ref{sec:blowup}. In the next Section~\ref{sec:smallk0} we shall have a more detailed look at approximation of slow manifolds and under which conditions one may be able to tackle the case $k_0\ra \infty$.

\section{Slow manifold approximation and Galerkin limits}
\label{sec:smallk0}

To understand the dynamical structure of the Galerkin ODEs, and in order to introduce the basic objects of multiple time scale dynamics and their notation, we analyze equation~\eqref{eq:k0_general} for $k_0=1,2$, omitting the higher order terms $H_k^u$ and $H_k^v$ for now.
\subsection{Truncation at first mode}
\label{sssec:truncate1}
Clearly, for $k_0 =1$, we have
\be
\label{eq:k01}
\begin{array}{lcl}
\partial_t u_1&=&  u_1^2 -  v_1^2 + 2a \varepsilon \mu,\\
\partial_t v_1&=& 2a \varepsilon  , \\
\partial_t \varepsilon &=& 0,
\end{array} 
\ee
such that we can describe the behaviour in the leading modes. Recall that the critical manifold~\cite{Fenichel4,Jones,KuehnBook} in standard form can be defined as the zero set of the fast vector field upon considering the limit $\varepsilon\ra 0$; we note that the terminology ``critical manifold'' is commonly used, even though more general algebraic varieties appear, which already provides a clear indication, why blowing-up can be useful~\cite{Hironaka1,Hironaka2}.
For~\eqref{eq:k01}, we focus on the branch $S_{\txta}^{-,1}$ of the critical manifold $S^1= \{u_1^2 = v_1^2\}$ given by
$$ S_{\txta}^{-,1} = \{ (u_1, v_1) \,:\,u_1 = v_1 < 0 \}.$$ 
The other three natural branches away from the origin $(u_1,v_1)=(0,0)$ can be treated analogously. A direct linearization of the fast vector field with respect to the fast variables yields the local stability of $S_{\txta}^{-,1}$. We obtain $\txtD_{u_1}(u_1^2-v_1^2+2a\varepsilon\mu)|_{S_{\txta}^{-,1}}=2u_1<0$, so $S_{\txta}^{-,1}$ is normally hyperbolic (as all eigenvalues of the fast subsystem linearization have non-zero real part) and also attracting (as all eigenvalues have negative real part). Fenichel's Theorem~\cite{Fenichel4,Jones,KuehnBook} provides the existence of a normally hyperbolic locally invariant slow manifold. The slow manifold can also be viewed as a center manifold of system~\eqref{eq:k01} due to its construction. We want to derive a parametrization $h_1(v_1, \varepsilon)$ of the center manifold at a given point on $S_{\txta}^-$, i.e.~at the equilibrium $\{u_1=v_1= c, \varepsilon=0\}$ for $c <0$ with $\left| c \right|$ sufficiently small. Understanding this parametrization will be important later on in our construction of the slow manifolds in the normally hyperbolic regime before the bifurcation point. We first change coordinates to shift the equilibrium to the origin
$$ \tilde u_1 = u_1 -c, \ \tilde v_1 = v_1 -c,$$
obtaining the equations
\be
\label{eq:k01_shifted}
\begin{array}{lcl}
\partial_t \tilde u_1&=&  (\tilde u_1+c)^2 -  (\tilde v_1+c)^2 + 2a \varepsilon \mu,\\
\partial_t \tilde v_1&=& 2a \varepsilon  , \\
\partial_t \varepsilon &=& 0.
\end{array} 
\ee
The Jacobian of the vector field for this ODE, at $\tilde u_1 = 0, \tilde v_1 = 0, \varepsilon =0$, reads
$$ A = \begin{pmatrix}
2c & - 2c & 2 a \mu \\
0 & 0& 2 a\\
0 & 0 & 0
\end{pmatrix},
$$ 
with eigenvalues $\lambda_1 = 2 c < 0$ and $\lambda_{2,3} =0$.
The Jordan normal form gives
\begin{equation*}
M^{-1} A M = \begin{pmatrix}
2c & 0 & 0 \\
0 & 0& 1\\
0 & 0 & 0
\end{pmatrix},
\end{equation*}
where
$$ M = \begin{pmatrix}
1 & 1 & \frac{1 - \mu}{2c} \\
0 & 1& 0\\
0 & 0 & \frac{1}{2a}
\end{pmatrix}, \ M^{-1} = \begin{pmatrix}
1 & - 1 & 2 a\left(\frac{\mu-1}{2c} \right) \mu \\
0 & 1& 0\\
0 & 0 & 2 a
\end{pmatrix}.
$$
The coordinate change
$$ M^{-1} \begin{pmatrix}
\tilde u_1 \\ \tilde v_1 \\ \varepsilon
\end{pmatrix}
= \begin{pmatrix}
y_1 \\ x_1 \\ x_2
\end{pmatrix}$$
gives 
\begin{align*}
\tilde u_1 &= y_1 + v_1 - 2 a \varepsilon \left( \frac{\mu -1}{2c} \right), \\
\tilde v_1 &= x_1,\\
\varepsilon &= \frac{x_2}{2a},
\end{align*}
such that, for $g = \mathcal{O}(2)$, we obtain
\begin{align}
y_1' &= 2c y_1 + g(x_1, x_2, y_1):= K_1(x_1,x_2,y_1), \\
\begin{pmatrix}
x_1' \\
x_2'
\end{pmatrix}
&= \begin{pmatrix}
0 & 1 \\
0 & 0
\end{pmatrix}
\begin{pmatrix}
x_1 \\
x_2
\end{pmatrix} =: \begin{pmatrix}
F_1(x_1, x_2) \\
F_2(x_1, x_2)
\end{pmatrix}.
\end{align}
This problem is now in standard form and we can approximate the center manifold via
$$ \hat h_1(x_1, x_2) = b_{11} x_1^2 + b_{12} x_1 x_2 + b_{22} x_2^2 + \mathcal{O}(3).$$
Solving the invariance equation 
$$ K_1 = F_1 \partial_{x_1} h_1 + F_2 \partial_{x_2} h_1,$$
gives the coefficients
\begin{equation} \label{eq:coeff_k01}
b_{11} =0, \ b_{12} = \frac{\mu -1}{2c^2}, \ b_{22} = \frac{(\mu-3)(\mu -1)}{8c^3}.
\end{equation}
Transforming back gives the parametrization
\begin{equation}\label{eq:cmtilde_k01}
\tilde u_1 = \tilde h_1 (\tilde v_1, \varepsilon) = \tilde v_1 - \frac{a(\mu -1)}{c} \varepsilon + \frac{a(\mu-1)}{c^2}\tilde v_1 \varepsilon - \frac{a^2(\mu-3)(\mu -1)}{2c^3} \varepsilon^2 + \mathcal{O}(3),
\end{equation}
and, hence,
\begin{equation}\label{eq:cm_k01}
 u_1 =  h_1 ( v_1, \varepsilon) =  v_1 - \frac{a(\mu -1)}{c} \varepsilon + \frac{a(\mu-1)}{c^2}(v_1-c) \varepsilon - \frac{a^2(\mu-3)(\mu -1)}{2c^3} \varepsilon^2 + \mathcal{O}(3),
\end{equation}
One can check by computation that, indeed, $\tilde h_1 (\tilde v_1, \epsilon)$ satisfies the invariance equation associated with~\eqref{eq:k01_shifted}, given by
\begin{equation} \label{eq:invar_k01_shifted}
 2 a \varepsilon \partial_{\tilde v_1} \tilde h_1 (\tilde v_1, \varepsilon) = (\tilde h_1 ( \tilde v_1, \varepsilon) + c)^2 - (\tilde v_1 +c)^2 + 2 a \varepsilon \mu + \mathcal{O}(3).
\end{equation}

\begin{figure}[ht]
        \centering
        \begin{subfigure}{0.5\textwidth}
        \centering
  		\begin{overpic}[width=0.9\linewidth]{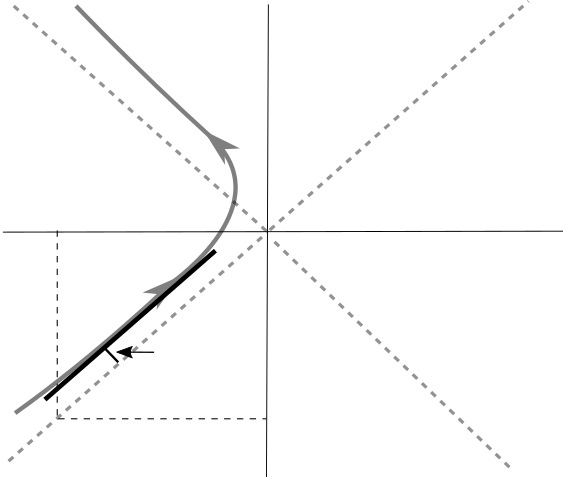}   
 \put(11,29){ $G_{\varepsilon, \zeta}$}
 \put(28,21){ $\mathcal O (\varepsilon)$}
 \put(90,39){ $u_1$}	
	\put(48,79){ $v_1$}
 \put(47,10){ $c$}
  \put(9,45){ $c$}
  \put(-4,6){ $S_{\txta}^{-,1}$}
        \end{overpic}
        \caption{$\mu < 1$}
        \label{k01smaller}
		\end{subfigure}%
        \begin{subfigure}{0.5\textwidth}
        \centering
  		\begin{overpic}[width=0.9\linewidth]{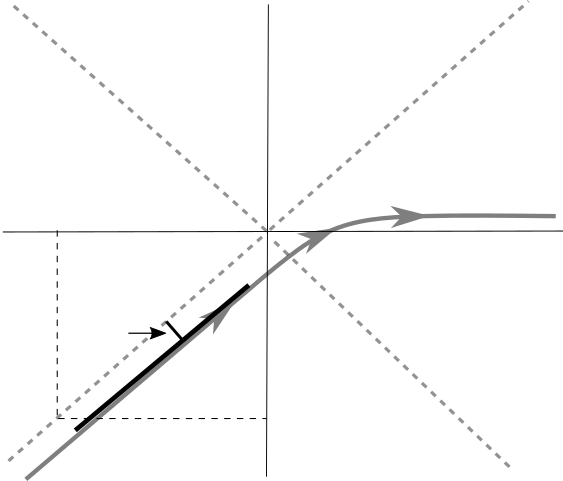}  
 \put(33,20){ $G_{\varepsilon, \zeta}$}
 \put(10,25){ $\mathcal O (\varepsilon)$}
 \put(90,39){ $u_1$}	
	\put(48,79){ $v_1$}
 \put(47,10){ $c$}
  \put(9,45){ $c$}
  \put(-5,6){ $S_{\txta}^{-,1}$}
        \end{overpic}
        \caption{ $\mu > 1$}
        \label{k01larger}
		\end{subfigure}
		\caption{ 
		    Depiction of the slow Galerkin manifolds close to  $S_{\txta}^{-,1}$  in the $(u_1, v_1)$-plane for $k_0 =1$, anchored at a fixed $u_1 = v_1 = c$ and parameterized by $u_1 = h_1(v_1, \varepsilon)$ \eqref{eq:cm_k01}, i.e.~$G_{\varepsilon, \zeta}$ for $ - \frac{9\pi^2}{4a} < \zeta^{-1} \omega_A < -\frac{\pi^2}{a}$ and $\varepsilon$ sufficiently small, when (a) $\mu < 1$ and (b) $\mu > 1$. The curves with arrows indicate the behaviour around the singularity at the origin as known from the classical planar ODE treatment, see e.g.~\cite{KruSzm4}, where the blow-up method is used for the extension of the slow manifolds.
      }
        \label{fig:k01}
\end{figure}

Hence, for this case, we have an explicit expression for the Galerkin parametrization $h^{\varepsilon, \zeta}_G (v_1) = h_1(v_1, \varepsilon) $, where $ - \frac{\pi^2}{a} < \zeta^{-1} \omega_A < - \frac{\pi^2}{4a}$ and $\varepsilon$ sufficiently small.
    In Figure~\ref{fig:k01}, we illustrate the first slow Galerkin manifolds for different values of $\mu$, also showing their extensions around the singularity at the origin, as described by the results in \cite{KruSzm4}. Note that the signs of the perturbations to first order in $\epsilon$ according to formula~\eqref{eq:cm_k01} (recalling $c<0$ and following $v_1 \geq c$) are consistent with the dynamically anticipated locations of the slow manifold: for $\mu < 1$ the dynamics tend to the left of the diagonal $\{u=v\}$, and for $\mu > 1$ the dynamics tend to the right of the diagonal $\{u=v\}$; see also \cite{EngelKuehn} for a detailed analysis of this behavior understood via time discretization. In fact, we oberve directly from formula~\eqref{eq:cm_k01} that for $\mu =1$ the dynamics can be predicted to stay on the diagonal  $\{u=v\}$; this is the canard case which we generally leave for future work in the situation of PDEs and Galerkin approximnations as it will require a detailed analysis in its own right.

To understand, how the Galerkin modes influence the dynamics, it is now instructive to also compute the next truncation level and compare results.

\subsection{Truncation at second mode}
\label{sssec:truncate2}
For $k_0 =2$, equation~\eqref{eq:k0_general} reads
\begin{equation} \label{eq:k02}
\begin{array}{lcl}
\partial_t u_1&=& u_1^2 - v_1^2 + 2a \varepsilon \mu +  u_2^2 - v_2^2 ,\\
\partial_t v_1&=& 2a\varepsilon , \\
\partial_t u_2&=& \hat \lambda_2 u_2 + 2 (u_2 u_1  -  v_2 v_1),\\
\partial_t v_2&=& \varepsilon \hat \lambda_2 v_2, \\
\partial_t \varepsilon &=& 0.
\end{array}
\end{equation} 
    We focus on the critical branch 
$$ S_{\txta}^{-,2} = \{ (u_1, v_1, u_2, v_2) \,:\,u_1 = v_1 < 0, \ u_2 = v_2 = 0\},$$
where we, again, fix $c < 0$ with $\left|c \right|$ small enough, and analyze the center manifold at the corresponding equilibrium.
Shifting 
to $\tilde u_1 = u_1 - c$ and $ \tilde v_1 = v_1 -c$ yields
\begin{equation} \label{eq:k02_shifted}
\begin{array}{lcl}
\partial_t \tilde u_1&=& (\tilde u_1+c)^2 - (\tilde v_1+c)^2 + 2a \varepsilon \mu +  u_2^2 - v_2^2 ,\\
\partial_t \tilde v_1&=& 2a\varepsilon , \\
\partial_t u_2&=& \hat \lambda_2 u_2 + 2 (u_2 (\tilde u_1 + c)  -  v_2 (\tilde v_1 + c)),\\
\partial_t v_2&=& \varepsilon \hat \lambda_2 v_2, \\
\partial_t \varepsilon &=& 0.
\end{array}
\end{equation}
We now consider the equilibrium at the origin, and obtain, compared to the case $k_0 =1$, the additional stable direction in $u_2$ with eigenvalue $\hat \lambda_2 + 2c < 0$ 
and center direction in $v_2$. Proceeding analogously to the case $k_0=1$ and transferring to standard form coordinates
\begin{align*}
\tilde u_1 &= y_1 + v_1 - 2 a \varepsilon \left( \frac{\mu -1}{2c} \right), \\
u_2 &= y_2 + \frac{2c}{2c + \hat \lambda_2}v_2,\\
\tilde v_1 &= x_1,\\
\varepsilon &= \frac{x_2}{2a},\\
v_2 &= x_3,
\end{align*}
we consider
\begin{align*}
\hat h_1(x_1, x_2, x_3) &= b_{11} x_1^2 + b_{12} x_1 x_2 + b_{22} x_2^2 + b_{13} x_1 x_3 + b_{23} x_2 x_3 + b_{33} x_3^2 + \mathcal{O}(3), \\
\hat h_2(x_1, x_2, x_3) &= c_{11} x_1^2 + c_{12} x_1 x_2 + c_{22} x_2^2 + c_{13} x_1 x_3 + c_{23} x_2 x_3 + c_{33} x_3^2 + \mathcal{O}(3).
\end{align*}
The associated invariance equation gives the coefficients
\begin{align} 
b_{11} &= b_{13}= b_{23} =0, \ b_{12} = \frac{\mu -1}{2c^2}, \ b_{22} = \frac{(\mu-3)(\mu -1)}{8c^3}, \ b_{33} = \frac{1}{2c} - \frac{2c}{(\hat \lambda_2 + 2c)^2},\label{eq:coeff_k02_b}\\
c_{11} &=  c_{22} = c_{33} = c_{12}= 0, \ c_{13} = \frac{2 \hat \lambda_2}{(\hat \lambda_2 + 2c)^2}, \ c_{23} = \frac{c \hat \lambda_2(2c+\hat \lambda_2)+4 ac(\mu-1)+2a \hat \lambda_2 \mu}{a(\hat \lambda_2 + 2c)^3}. \label{eq:coeff_k02_c}
\end{align}
Transforming back gives the center manifold approximation in form of the maps
\begin{align}\label{eq:cmtilde_k02_1}
\tilde u_1 &= \tilde h_1 (\tilde v_1, \varepsilon, v_2) \nonumber\\
&= \tilde v_1 - \frac{a(\mu -1)}{c} \varepsilon + \frac{a(\mu-1)}{c^2}\tilde v_1 \varepsilon - \frac{a^2(\mu-3)(\mu -1)}{2c^3} \varepsilon^2 + \left( \frac{1}{2c} - \frac{2c}{(\hat \lambda_2 + 2c)^2}\right) v_2^2 + \mathcal{O}(3),
\end{align}
and
\begin{align}\label{eq:cmtilde_k02_2}
\tilde u_2 &= \tilde h_2 (\tilde v_1, \varepsilon, v_2) \nonumber \\
&= \frac{2c}{2c + \hat \lambda_2} v_2 + \frac{2 \hat \lambda_2}{(\hat \lambda_2 + 2c)^2} \tilde v_1 v_2 + 2 \frac{c \hat \lambda_2(2c+\hat \lambda_2)+4 ac(\mu-1)+2a \hat \lambda_2 \mu}{(\hat \lambda_2 + 2c)^3} v_2 \varepsilon +  \mathcal{O}(3).
\end{align}
One can check by computation that, indeed, $\tilde h= (\tilde h_1, \tilde h_2)$ satisfies the system of invariance equations associated with~\eqref{eq:k02_shifted}, given by
\begin{equation} \label{eq:invar_k02_shifted}
\begin{array}{lcl}
 2 a \varepsilon \partial_{\tilde v_1} \tilde h_1 + \varepsilon \hat \lambda_2 v_2 \partial_{\tilde v_2} \tilde h_1  &=& (\tilde h_1 + c)^2 - (\tilde v_1 +c)^2 + 2 a \varepsilon \mu + (\tilde h_2)^2 - v_2^2 + \mathcal{O}(3),\\
 2 a \varepsilon \partial_{\tilde v_1} \tilde h_2 + \varepsilon \hat \lambda_2 v_2 \partial_{\tilde v_2} \tilde h_2 &=& \hat \lambda_2 \tilde h_2 + 2 (\tilde h_2 (\tilde h_1 + c)  -  v_2 (\tilde v_1 + c)) + \mathcal{O}(3).
 \end{array}
\end{equation}

\begin{figure}[ht]
        \centering
        \begin{subfigure}{0.5\textwidth}
        \centering
  		\begin{overpic}[width=0.95\linewidth]{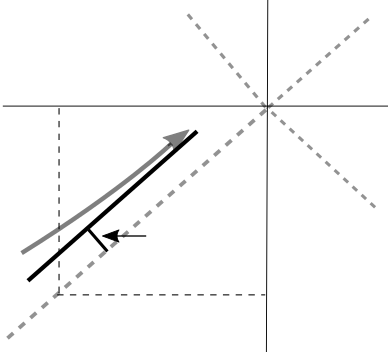}   
        \put(92,60){ $u_1$}	
	\put(69,85){ $v_1$}
 \put(48,58){ $G_{\varepsilon, \zeta}$}
 \put(39,29){\small $\mathcal O \left(\varepsilon - v_2^2\right)$}
 \put(68,14){ $c$}
  \put(13,65){ $c$}
  \put(5,4){ $S_{\txta}^{-,2}$}
        \end{overpic}
        \caption{$\mu < 1$}
        \label{k02smaller}
		\end{subfigure}%
        \begin{subfigure}{0.5\textwidth}
        \centering
  		\begin{overpic}[width=0.95\linewidth]{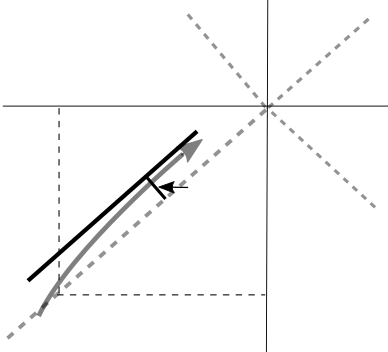}  
\put(92,60){ $u_1$}	
	\put(69,85){ $v_1$}
 \put(48,58){ $G_{\varepsilon, \zeta}$}
 \put(49,41){\small $\mathcal O \left(\varepsilon - v_2^2\right)$}
 \put(68,14){ $c$}
  \put(13,65){ $c$}
  \put(5,3){ $S_{\txta}^{-,2}$}
        \end{overpic}
        \caption{ $\mu > 1$}
        \label{k02larger}
		\end{subfigure}
		\caption{ Depiction of the slow Galerkin manifolds close to  $S_{\txta}^{-,2}$  projected to the $(u_1, v_1)$-plane for $k_0 =2$, anchored at a fixed $u_1 = v_1 = c, u_2 = v_2 =0$ and parameterized by $u_1 = h_1(v_1, v_2, \varepsilon)$ \eqref{eq:cm_k02_1}, i.e.~$G_{\varepsilon, \zeta}$ for $ - \frac{9\pi^2}{4a} < \zeta^{-1} \omega_A < -\frac{\pi^2}{a}$  and $\varepsilon$ sufficiently small, when (a) $\mu < 1$ and (b) $\mu > 1$. 
      The curves with arrows indicate the flow close to the manifolds. We note that the perturbation depending on $v_2^2$ tends to the left and, hence, may lead to $G_{\varepsilon, \zeta}$ lying to the left of $S_{\txta}^{-,2}$ also for $\mu > 1$ (b), depending on the sizes of $\varepsilon$ and $v_2^2$.}
        \label{fig:k02}
\end{figure}

The relations $u_1 = h_1 (v_1,  v_2, \varepsilon)$ and $u_2 = h_2 ( v_1, v_2, \varepsilon)$ are obtained by transforming back into original coordinates, 
    yielding
    \begin{align}\label{eq:cm_k02_1}
 u_1 &=  h_1 ( v_1, \varepsilon, v_2) \nonumber\\
&=  v_1 - \frac{a(\mu -1)}{c} \varepsilon + \frac{a(\mu-1)}{c^2} (v_1 - c) \varepsilon - \frac{a^2(\mu-3)(\mu -1)}{2c^3} \varepsilon^2 + \left( \frac{1}{2c} - \frac{2c}{(\hat \lambda_2 + 2c)^2}\right) v_2^2 + \mathcal{O}(3),
\end{align}
and
\begin{align}\label{eq:cm_k02_2}
 u_2 &=  h_2 ( v_1, \varepsilon, v_2) \nonumber \\
&= \frac{2c}{2c + \hat \lambda_2} v_2 + \frac{2 \hat \lambda_2}{(\hat \lambda_2 + 2c)^2} (v_1-c) v_2 + 2 \frac{c \hat \lambda_2(2c+\hat \lambda_2)+4 ac(\mu-1)+2a \hat \lambda_2 \mu}{(\hat \lambda_2 + 2c)^3} v_2 \varepsilon +  \mathcal{O}(3).
\end{align}
Hence, also for this case, we have an explicit expression for the Galerkin parametrization $h^{\varepsilon, \zeta}_G(v_1,v_2) = \left( h_1(v_1,  v_2, \varepsilon), h_2(v_1,  v_2, \varepsilon) \right) $, where $ - \frac{9\pi^2}{4a} < \zeta^{-1} \omega_A < -\frac{\pi^2}{a}$ and $\varepsilon$ sufficiently small.
    For $u_1$, note the additional perturbation term depending on $v_2^2$ which has a negative sign for any $\left| c \right|$ small enough. In Figure~\ref{fig:k02}, we illustrate the effect of $v_2 \neq 0$ for the cases $\mu <1$ , where it is qualitatively negligible since the position of the slow manifold does not change with respect to the diagonal, and for $\mu > 1$, where it has the potential effect of shifting $G_{\varepsilon, \zeta}$ to the left of the diagonal for sufficiently small $\varepsilon > 0$. This fact gives a first hint at the subtleties we have to deal with when analyzing the continuation around the origin in the situation of higher Galerkin modes being at play. 
    As we loose the hyperbolic perturbation procedure when $c \to 0$, we will need a careful application of the blow-up method in order to determine the behaviour when passing the singularity at the origin, bounding also sufficiently small initial conditions for additional modes like $v_2$, as conducted in Section~\ref{sec:blowup}.
    Note that, additionally, we now also have to parameterize $u_2 \approx 0$ close to $v_2 \approx 0$, $v_1 \approx c$, yielding the control on keeping $u_2$ close to $0$ in terms of the slow modes.

    In the situation of $\mu=1$, which we do not discuss here in detail due to the canard phenomenon, the dependence on $v_2$ has an interesting qualitative effect, apparently breaking the symmetry of the canard: as already stated above, we leave a thorough analysis of this situation for future work.

The calculation for the case $k_0=3$ can be found in Appendix~\ref{ap:truncate3}. 
Higher-order approximations can be calculated similarly, as shown in the next subsection.

\subsection{Truncation at arbitrary mode}
We now consider the general case $k_0 \geq 1$, again skipping the higher-order terms $H_k^u$ and $H_k^v$, which do not have a significant impact on the following approximation of the center manifolds up to second order. Focusing on the critical branch 
$$ S_{\txta}^{-,k_0} = \{ (u_1, v_1, u_2, v_2, \dots, u_{k_0}, v_{k_0}) \,:\,u_1 = v_1 < 0, \ u_k = v_k = 0, \ 2\leq k \leq k_0 \},$$
we pick $c < 0$ with $\left|c \right|$ small enough, yet still a fixed constant, and analyze the center manifold at the corresponding equilibrium.
Note the accordance with $G_{\varepsilon, \zeta}$ for appropriate $\zeta$.

In more details, we make the change of variables 
$$ \tilde u_1 = u_1 - c, \ \tilde v_1 = v_1 -c,$$
such that, skipping the tilde notation in the following, we obtain the system of equations, with equilibrium at the origin,
\begin{subequations}
\label{eq:k0_general_shift}
\begin{align}
\partial_t  u_1&= ( u_1+c)^2 - ( v_1+c)^2 + 2a \varepsilon \mu +  \sum_{j=2}^{k_0} \left(u_j^2 - v_j^2 \right),\\
\partial_t v_1&= 2a\varepsilon , \\
\partial_t u_k&= \hat \lambda_k u_k + 2 (u_k ( u_1 +c)  -  v_k ( v_1 +c)) + \sqrt{2}\sum_{i,j=2}^{k_0} \alpha_{i,j}^k (u_i u_j - v_i v_j), \ 2 \leq k\leq k_0,\\
\partial_t v_k&= \varepsilon \hat \lambda_k v_k, \ 2 \leq k \leq k_0, \\
\partial_t \varepsilon &= 0.
\end{align} 
\end{subequations}
For $c<0$, the origin is a center-stable equilibrium (center-unstable for $c>0$) with center manifold  
\begin{equation} \label{center_stable_before_blowup_k0}
\hat M_c^{k_0} = \left\{\left( u^{k_0}, ( v^{k_0}, \varepsilon)\right) \in  \mathbb{R}^{k_0} \times \mathbb{R}^{k_0+1}\,:\, \left( u_1^{k_0}, \dots, u_{k_0}^{k_0} \right) =  h_c^{k_0}( v^{k_0}, \varepsilon) \right\}.
\end{equation}
For a neighbourhood $\mathcal{U}(0) \subset L^2([-a,a])$, we can view the parametrization map $ h^{k_0}$ as a map
$$  h_c^{k_0} : \mathcal{U}(0) \times [0, \varepsilon_1] \to L^2([-a,a]), \ ( v, \varepsilon) \mapsto u^{k_0} $$
for some fixed $\varepsilon_1 > 0$, via identifying
\begin{align*}
 v &= \sum_{k=1}^{\infty} v_k e_k, \\
 u^{k_0}  &= \sum_{k=1}^{k_0} u_k^{k_0} e_k,
\end{align*}
and denoting
$ v^{k_0} = \left( v_1, \dots,v_{k_0} \right)$.
Hence, we can consider $\hat M_c^{k_0}$ also as a Banach manifold
\begin{equation} \label{center_stable_before_blowup_k0_L2}
 M_c^{k_0} = \left\{\left( u, ( v, \varepsilon)\right) \in  L^2([-a,a]) \times \left(\mathcal{U}(0) \times [0, \varepsilon_1]\right)\,:\, u =  h_c^{k_0}( v, \varepsilon) \right\}.
\end{equation}
Replacing $h_c^{k_0}$ by $h^{k_0}$ in the following for ease of notation,  the invariance equation for system~\eqref{eq:k0_general_shift} and $h^{k_0} = (h_1^{k_0}, \dots, h_{k_0}^{k_0})$ reads
\begin{subequations}\label{eq:invar_k0_general}
\begin{align}
 2 a \varepsilon \partial_{ v_1}  h_1^{k_0}  + \varepsilon \sum_{i=2}^{k_0} \hat \lambda_i v_i \partial_{v_i}  h_1^{k_0}    &= (h_1^{k_0}  + c)^2 - (v_1 +c)^2 + 2 a \varepsilon \mu + \sum_{i=2}^{k_0} \left((h_i^{k_0})^2 - v_i^2\right),\label{eq:invar_k0_general_1} \\
 2 a \varepsilon \partial_{ v_1}  h_k^{k_0}  + \varepsilon \sum_{i=2}^{k_0} \hat \lambda_i v_i \partial_{v_i}  h_k^{k_0} &= \hat \lambda_k  h_k^{k_0} + 2 ( h_k^{k_0} ( h_1^{k_0} + c)  -  v_2 ( v_1 + c)) \nonumber\\
 &+\sqrt{2}\sum_{i,j=2}^{k_0} \alpha_{i,j}^k \left( h_i^{k_0} h_j^{k_0} - v_i v_j \right), \ \text{ for all } 2 \leq k \leq k_0. \label{eq:invar_k0_general_k}
\end{align}
\end{subequations}
First of all, we show the following lemma, generalizing the observations from Section~\ref{sec:smallk0}:
\begin{lemma} \label{lem:solution_of_invariance_equation}
The invariance equation~\eqref{eq:invar_k0_general} has the solution (under changing back to coordinates $u_1 = \tilde u_1 +c, v_1 = \tilde v_1 +c$)
\begin{align}
u_1^{k_0} &= h_1^{k_0}( v^{k_0}, \varepsilon) = v_1 - \frac{a(\mu -1)}{c} \varepsilon + \frac{a(\mu-1)}{c^2} (v_1 - c) \varepsilon \nonumber\\
&- \frac{a^2(\mu-3)(\mu -1)}{2c^3} \varepsilon^2 + \sum_{j=2}^{k_0} \left( \frac{1}{2c} - \frac{2c}{(\hat \lambda_j + 2c)^2}\right) v_j^2 + \mathcal{O}(3), \label{eq:cmtilde_k0gen_1} \\
u_k^{k_0} &= h_k^{k_0}( v^{k_0}, \varepsilon) = \frac{2c}{2c + \hat \lambda_k} v_k +  2 \frac{c \hat \lambda_k(2c+\hat \lambda_k)+4 ac(\mu-1)+2a \hat \lambda_k \mu}{(\hat \lambda_k + 2c)^3} v_k \varepsilon \nonumber \\
&+ \frac{2 \hat \lambda_k}{(\hat \lambda_k + 2c)^2} (v_1 -c) v_k + \frac{C_k}{(2c+\hat \lambda_k)} \sum_{i,j=2}^{k_0} \beta_{i,j}^k v_i v_j +  \mathcal{O}(3), \label{eq:cmtilde_k0gen_k}
\end{align}
where $ \left| C_k \right| < C$ for a constant $C \in \mathbb{R}$, and $ \left| \beta_{i,j}^k \right| \in [0,1]$ with 
$$\beta_{i,j}^k \neq 0 \ \text{ if and only if } \ i+j-k =1 \lor k-\left|i-j\right| =1.$$
\end{lemma}
\begin{proof}
For $k=1$ this is just a straight-forward calculation by inserting formula~\eqref{eq:cmtilde_k0gen_1} and 
$$\left(h_i^{k_0}\right)^2 = \left(\frac{2c}{2c + \hat \lambda_i} v_i\right)^2 + \mathcal{O}(3)$$
into equation~\eqref{eq:invar_k0_general_1}.

For $k \geq 2$, this is also a direct calculation, inserting formulas~\eqref{eq:cmtilde_k0gen_1} and~\eqref{eq:cmtilde_k0gen_k} into equation~\eqref{eq:invar_k0_general_k}. Note that the factor $(2 c + \hat \lambda_k)^{-1}$ in each term, is directly inherited from the terms $\hat \lambda_k  h_k^{k_0}$ and $ 2 h_k^{k_0} c$ on the right hand side of equation~\eqref{eq:invar_k0_general_k}. The fact that we can find a general constant $C > \left| C_k \right|$ can easily be seen from the fact that all other coefficients in equation~\eqref{eq:cmtilde_k0gen_k} are uniformly bounded and $\alpha_{i,j}^k \in [0,1]$.
\end{proof}
In other words, we can give a general explicit expression for the Galerkin parametrization up to second order, given by 
$$h^{\varepsilon, \zeta}_G(v^{k_0}) = \left( h_1^{k_0}(v^{k_0}, \varepsilon), h_2^{k_0}(v^{k_0}, \varepsilon), \dots, h_{k_0}^{k_0}(v^{k_0}, \varepsilon),  \right),$$
where
 \begin{equation*}
	- \frac{{(k_0 + 1)}^2 \pi^2}{4a^2} < \zeta^{-1} \omega_A \le - \frac{k_0^2 \pi^2}{4 a^2}, 
\end{equation*}
and $\varepsilon$ sufficiently small.

We briefly summarize some important observations from these explicit computations of the center/slow manifolds in the normally hyperbolic regime before the fast subsystem bifurcation point at $(u,v)=(0,0)$.
Comparing the different Galerkin orders, we observe that the basic structure from the first modes persists throughout the different orders. 
    Formulas~\eqref{eq:cmtilde_k0gen_1} and~\eqref{eq:cmtilde_k0gen_k} exhibit the same qualitative phenomena as explained for the case $k_0 = 2$; terms of $\mathcal{O} \left( \sum_{j=2}^{k_0} v_j^2\right)$ perturb the positioning of the manifolds $G_{\varepsilon, \zeta}$ projected to the $(u_1, v_1))$-plane (cf.~Figure~\ref{fig:k02}), where $k_0$ corresponds with $\zeta$ as above, and the size of $u_k \approx 0$, $k \geq 2$, is controlled by the size of the slower modes $v_k \approx 0$.
 When $c\nearrow 0$, the expansions diverge, as expected, since normal hyperbolicity breaks down at $c=0$. In fact, this breakdown necessitates the use of a nonlinear method, such as the blow-up method, to deal with the invariant manifolds passing near the singularity. 
    The blow-up analysis in Section~\ref{sec:blowup}, in particular the analysis of the rescaling chart in Section~\ref{sec:chart2}, will explicitly bound perturbation effects of the slow higher modes $v_2, \dots, v_{k_0}$ and keep the fast higher modes $u_2, \dots, u_{k_0}$ sufficiently small, extending the slow Galerkin manifolds parametrized in Lemma~\ref{lem:solution_of_invariance_equation} around the singularity at the origin.
 
Furthermore, we note that many additional terms scale with $(\hat\lambda_{k})^{-1}$ in the limit when $k\ra \I$. Since we have full control over the eigenvalues, it is easy to see that $(\hat\lambda_{k})^{-1}\ra 0$ as $k\ra \I$. This implies that we can conjecture the existence of a limiting manifold as the Galerkin truncation level is increased as long as the $\mathcal{O}(3)$ terms are ignored; we remark that this argument seems to hold on even far more general domains as long as we can rely on Weyl's law for the eigenvalues of the Laplacian.

\subsection{Convergence to infinite-dimensional manifold}

Taking $\mathcal{U}(0) \subset L^2([-a,a])$ to be a small neighbourhood, we can now show the following proposition, which naturally also holds when transforming back via $ u_1 = \tilde u_1 + c, \ v_1 = \tilde v_1 + c$ into a neighbourhood of $\mathcal{U}(v^*) \subset L^2([-a,a])$ for $v^* \equiv c$:
\begin{proposition} \label{prop:convergence_beforeblowup}
For any $v \in \mathcal{U}(0) \subset L^2([-a,a])$, we have that
\begin{align}
\left|u_1^{k_0}\right| &= \left|h_1^{k_0}( v^{k_0}, \varepsilon) \right| \leq C_1(c) \left(\|v\|_2^2 + (1+\varepsilon)^2\right) + \mathcal{O}(3), \label{eq:bound_u1}\\
\left|u_k^{k_0}\right| &= \left|h_k^{k_0}( v^{k_0}, \varepsilon) \right| \leq \frac{C_2(c)}{\left|\hat \lambda_k\right|} (\|v\|_2^2 + (1+\varepsilon)) + \mathcal{O}(3). \label{eq:bound_uk}
\end{align}
for all $k_0 \in \mathbb{N}$, where $C_1(c), C_2(c) > 0$ are constants.
\end{proposition}
\begin{proof}
Using the fact that 
$$\sup_{v \in \mathcal{U}(0), \,  k \in \mathbb N} \left|v_k \right| < \infty, $$ 
we only need to take care of the terms
$$\sum_{j=2}^{k_0} \left( \frac{1}{2c} - \frac{2c}{(\hat \lambda_j + 2c)^2}\right) v_j^2$$
in equation~\eqref{eq:cmtilde_k0gen_1} and
$$ \frac{C_k}{(2c+\hat \lambda_k)} \sum_{i,j=2}^{k_0} \beta_{i,j}^k v_i v_j$$
in equation~\eqref{eq:cmtilde_k0gen_k}, to obtain the bounds~\eqref{eq:bound_u1} and~\eqref{eq:bound_uk}.
Since we have
$$ \|v\|_2^2 = \sum_{j=1}^{\infty} \left| v_j \right|^2, $$
we immediately get~\eqref{eq:bound_u1}. Furthermore, recall that $\left| C_k \right| < C < \infty$ and 
$ \beta_{i,j}^k \in \{0,1\}$ with $\beta_{i,j}^k =1$ if and only if $i+j-k =1 \lor k-\left|i-j\right| =1$.
Hence, we obtain by using H{\"{o}}lder's inequality
\begin{align*}
&\left| \frac{C_k}{\left|2c+\hat \lambda_k\right|} \sum_{i,j=2}^{k_0} \beta_{i,j}^k v_i v_j \right| \\
&\leq \frac{C}{\left|2c+\hat \lambda_k\right|} \left( \sum_{i=2}^{k-1} \left| v_i v_{k+1-i} \right| + 2 \sum_{i=k+1}^{k_0} \left| v_i v_{i+1-k} \right| \right) \\
&\leq \frac{C}{\left|2c+\hat \lambda_k\right|} \left( \left(\sum_{i=2}^{k-1} \left| v_i \right|^2 \right)^{1/2} \left(\sum_{i=2}^{k-1} \left| v_{k+1-i}  \right|^2 \right)^{1/2} + 2 \left(\sum_{i=k+1}^{k_0} \left| v_i\right|^2\right)^{1/2} \left(\sum_{i=k+1}^{k_0} \left| v_{i+1-k}\right|^2\right)^{1/2}  \right)\\
&\leq 3 \frac{C}{\left|2c+\hat \lambda_k\right|}  \|v\|_2^2.
\end{align*}
This gives~\eqref{eq:bound_uk}. 
\end{proof}
The following remark discusses the consequences of this observation if we ignore higher-order effects of $\varepsilon$. Note that for fixed $k_0$, this is always justifiable for $\varepsilon$ sufficiently small. However, when $k_0 \to \infty$, we also need to consider $\varepsilon \to 0$, such that the following limits are not necessarily valid in a rigorous sense. At this stage, we are not able to close this gap but leave this remark as a potentially helpful insight for simulations as well as theoretical progress.
\begin{remark}
\label{rem:convergence_beforeblowup}
 \begin{enumerate}
     \item[(i)] If we ignore the $\mathcal{O}(3)$ terms (which is possible for each fixed $k_0$ but potentially problematic in the limit), we may immediately deduce that there exists $\varepsilon_1>0$ such that for all $\varepsilon\in[0,\varepsilon_1]$,
 the sequence $(u^n)_{n \in \mathbb{N}} \subset L^2([-a,a])$ is bounded and contains a weakly convergent subsequence. This follows by observing with~\eqref{eq:bound_uk} that
$$ \|u^n\|_2^2 = \sum_{k=1}^{n} \left|u_k^n\right|^2 \leq K_1 + K_2 \sum_{k=1}^{\infty} \left|\hat \lambda_k\right|^{-2} < \infty,$$
where $K_1, K_2 > 0$ are constants for fixed $v$. Thus, the existence of a weakly convergent subsequence follows with the Banach-Eberlein-Smulian theorem \cite[Theorem 5.14-4]{Ciarlet2013}.

\item[(ii)] Even more strongly, we then have $h^{n}( v, \varepsilon) = u^{n} \to u^*=:h( v, \varepsilon) $ as $n \to \infty$ in $L^2([-a,a])$, uniformly over all $(v, \varepsilon) \in \mathcal{U}(0) \times [0, \varepsilon_1]$. 
This can be seen by showing that $(u^n)_{n \in \mathbb{N}}$ is a Cauchy sequence in $L^2([-a,a])$:
For that, let $\delta > 0$. Then it follows immediately from~\eqref{eq:bound_uk} that there is an $N_1 \in \mathbb{N}$ such that for all $N \geq N_1$, $m, n \in \mathbb{N}$ and a constant $K >0$ depending on $\|v\|$,
$$ \sum_{k=N}^{\infty} \left|u_k^n - u_k^m \right|^2 \leq \sum_{k=N}^{\infty} \left(\left|u_k^n\right| + \left| u_k^m \right| \right)^2 \leq K \sum_{k=N}^{\infty} \left| \hat \lambda_k \right|^{-2} < \frac{\delta}{3}. $$
Furthermore, we deduce from $v \in L^2([-a,a])$ and formula~\eqref{eq:bound_u1} that there is an $N_2 \in \mathbb{N}$ such that for all $n \geq m \geq N_2 > N_1$
$$ \left|u_1^n - u_1^m \right|^2 = \left| \sum_{j=m}^{n} \left( \frac{1}{2c} - \frac{2c}{(\hat \lambda_j + 2c)^2}\right) v_j^2 \right|^2 < \frac{\delta}{3}.$$
Additionally, we observe for $n > m > k\geq 2$ that
\begin{align*}
\left|u_k^n - u_k^m \right| &= \left| \frac{C_k}{(2c+\hat \lambda_k)} \right| \left| \sum_{i,j=2}^{n} \beta_{i,j}^k v_i v_j -\sum_{i,j=2}^{n} \beta_{i,j}^k v_i v_j \right| \\
&\leq  \left| \frac{C}{(2c+\hat \lambda_k)} \right| 2 \sum_{i=m+1}^{n} \left| v_i v_{i-k+1} \right| \\
&\leq   \left| \frac{2 C}{(2c+\hat \lambda_k)} \right| \|v\|_2 \left(\sum_{i=m+1}^{n} \left| v_i \right|^2 \right)^{1/2}.
\end{align*} 
Hence, there is an $N^* \geq N_2$ such that for $n \geq m\geq N^*$
$$ \sum_{k=2}^{N_1} \left|u_k^n - u_k^m \right|^2 \leq (N_1-1)\left| \frac{2 C}{(2c+\hat \lambda_2)} \right|^2 \|v\|_2^2 \left(\sum_{i=m+1}^{n} \left| v_i \right|^2 \right) < \frac{\delta}{3}.$$
Hence, we obtain that for $n \geq m\geq N^*$
$$ \|u^n - u^m \|_2^2 = \sum_{k=1}^{\infty} \left|u_k^n - u_k^m \right|^2 = \left|u_1^n - u_1^m \right|^2 + \sum_{k=2}^{N_1} \left|u_k^n - u_k^m \right|^2 + \sum_{k=N_1 + 1}^{\infty} \left|u_k^n - u_k^m \right|^2 <  \frac{\delta}{3} +  \frac{\delta}{3} +  \frac{\delta}{3} = \delta.$$
Since every estimate only depends on $\|v\| \in \mathcal{U}(0)$ which has a uniform bound, the convergence is uniform.
\end{enumerate}
\end{remark}

\begin{remark}
\label{rem:Mcconvergence}
If we are still in the situation of Remark~\ref{rem:convergence_beforeblowup}, i.e.~ignoring the $\mathcal{O}(3)$ terms, and take any given $c \in \mathbb R$, we may write $h_c( v, \varepsilon)$ for the limit in Remark~\ref{rem:convergence_beforeblowup} (ii) and define the infinite-dimensional Banach manifold 
\begin{equation} \label{center_stable_before_blowup}
M_c := \left\{\left( u, ( v, \varepsilon)\right) \in  L^2([-a,a]) \times \left( \mathcal{U}(0) \times [0, \varepsilon_1] \right)\,:\, u =  h_c( v, \varepsilon) \right\}.
\end{equation}
Introducing the Hausdorff distance $d_{\textnormal{H}}$ on sets $A, B \subset  L^2([-a,a]) \times (L^2([-a,a]) \times [0, \varepsilon_1]) $ via
$$d_{\textnormal{H}}(A,B) = \max \left\{ \sup_{x \in A} \inf_{y \in B} \|x-y\|_{L^2}, \ \sup_{y \in B} \inf_{x \in A} \|x-y\|_{L^2}  \right\}, $$
we may then deduce the following, where the manifolds $M_c^n$ may be slightly shifted due to the higher order terms $H_k^u, H_k^v$:
%
Denoting by $S_{\txta, \varepsilon}^{-,k_0}$ the lower left branch ($u_1<0$,$v_1<0$) of the slow Fenichel manifold for the fast-slow system of Galerkin ODEs~\eqref{eq:k0_general}, we have:
\begin{enumerate}
\item The manifolds $M_c^n$ \eqref{center_stable_before_blowup_k0_L2} converge to $M_c$~\eqref{center_stable_before_blowup} as $n \to \infty$ in Hausdorff distance, i.e.
$$ d_{\textnormal{H}}(M_c^n,M_c) \to 0.$$
\item For any $\varepsilon^* \in [0, \varepsilon_1]$, the manifolds $M_c^n \cap \{\varepsilon=\varepsilon^*\}$ and $S_{\txta, \varepsilon^*}^{-,n}$ coincide locally and, in particular, we have convergence of $S_{\txta, \varepsilon^*}^{-,n}$ to a Banach manifold $S_{\txta, \varepsilon^*}^{-}$ as $n \to \infty$ in Hausdorff distance.
\end{enumerate}
The first part is straightforwardly following from Remark~\ref{rem:convergence_beforeblowup} due to the uniform convergence of $(h_c^n(v, \varepsilon), v, \varepsilon) \to (h_c(v, \varepsilon), v, \varepsilon)$. The second part follows from classical Fenichel theory by identifying the center manifold $M_c^n \cap \{\varepsilon=\varepsilon^*\}$ and the slow manifold $S_{\txta, \varepsilon^*}^{-,n}$ as usual~\cite{Fenichel4}. Then we simply apply the convergence result for $M_c^n$.

The manifolds $M_c^n$, and the limiting object $M_c$, are then naturally transformed back into a neighbourhood $\mathcal{U}(u^*)  \times \mathcal{U}(v^*) \subset L^2([-a,a]) \times L^2([-a,a])$ for $u^* = v^* \equiv c$  via $ u_1 = \tilde u_1 + c, \ v_1 = \tilde v_1 + c$. Analogous statements
can be made around the branches
$$ S_{\txta}^{+,k_0} = \{ (u_1, v_1, u_2, v_2, \dots, u_{k_0}, v_{k_0}) \,:\,u_1 = - v_1 < 0, \ u_k = v_k = 0, \ 2\leq k \leq k_0 \},$$
for the corresponding objects $S_{\txta, \varepsilon}^{+,k_0}$, $S_{\txta, \varepsilon}^{+}$ and, for $c < 0$ with $\left|c \right|$ small, the center manifolds $\tilde M_c^n$ and $\tilde M_c$ corresponding with $u_1 = - v_1 = c$.
\end{remark}

There are two important aspects concerning the dynamical interpretation of the manifold $S_{\txta, \varepsilon^*}^{-}$ in Remark~\ref{rem:Mcconvergence}, when taking $\varepsilon^* \in [0, \varepsilon_1]$, i.e.~for some fixed upper bound $\varepsilon_1$. First, the observations in Remark~\ref{rem:convergence_beforeblowup} only hold up to terms of order $\mathcal{O}(3)$, i.e., we cover the dynamics up to a second-order approximation. 
These terms potentially grow as $k_0 \to \infty$ such that, in fact, the uniform bounds vanish and we observe $\varepsilon_1 = \varepsilon_1(k_0)\ra 0$ as $k_0\ra \I$.
Second, recall that the notion of a fast-slow splitting disappears for $k_0 \to \infty$ and fixed $\varepsilon$ since $\left|\lambda_{k} \right| \to \infty$. 
Hence, in full rigour, we have to work with the objects from Section~\ref{sec:manifolds}, i.e.~the $\zeta$ and by that $k_0$ parametrized slow manifolds $S_{\varepsilon, \zeta}$ and its Galerkin approximations
$G_{\varepsilon, \zeta}$, see Proposition~\ref{Prop:slow_Galerkin}.
Hence, in the following we track the dynamics of $S_{\varepsilon, \zeta}$ by a full analysis of the dynamics of the corresponding Galerkin manifold $G_{\varepsilon, \zeta}$ near the origin. 


\section{Blow-up analysis through the singularity}
\label{sec:blowup}
Consider again system~\eqref{eq:k0_general} with equilibrium at $u_1 = v_1 =c$, $u_k=v_k= \varepsilon=0$, $k \geq 2$, $c < 0$. The equilibrium looses hyperbolicity at $c=0$, i.e.~at the origin, implying that, in addition to the slow variables, also the first mode $u_1$ gives a center direction at this point. This is a situation, where geometric desingularization via the blow-up method has proven to be extremely useful; see~\cite{Du93,JardonKuehn,KuehnBook} for detailed introductions and surveys of blow-up methods for general ODEs and their particular effectiveness for multiple time scale ODEs.

To apply the blow-up to our system~\eqref{eq:k0_general}, it is first helpful to work with a scaled version of the eigenvalues for better grasping the role of $a > 0$. We write
\begin{equation} \label{eq:express_in_a}
\hat \lambda_k = - \pi^2 (k-1)^2 2^{-3/2} a^{-3/2} =: b_k a^{-3/2}.
\end{equation} 
For the following, keep in mind that $b_k$ is monotonically decreasing with
 $$b_k \to  -\infty \ \text{ and } \ \sum_{k=1}^{\infty} \left|b_k^{-1}\right| < \infty.$$
We write system~\eqref{eq:k0_general}, for any $k_0 \in \mathbb{N}$, as
\begin{subequations} \label{eq:k0_general_forblowup}
\begin{align}
\partial_t u_1&= u_1^2 - v_1^2 + 2a \varepsilon \mu +  \sum_{j=2}^{k_0} \left(u_j^2 - v_j^2 \right) + H_1^u,\\
\partial_t v_1&= 2a\varepsilon + H_1^v, \\
\partial_t u_k&= b_k a^{-3/2} u_k + 2 (u_k u_1  -  v_k v_1) + \sqrt{2}\sum_{i,j=2}^{k_0} \alpha_{i,j}^k (u_i u_j - v_i v_j) + H_k^u, \ 2 \leq k\leq k_0,\\
\partial_t v_k&= \varepsilon b_k a^{-3/2} v_k + H_k^v, \ 2 \leq k \leq k_0, \\
\partial_t \varepsilon &= 0,\\
\partial_t a &= 0,
\end{align} 
\end{subequations}
where $\alpha_{i,j}^k \in [0, 1]$ with $\alpha_{i,j}^k \neq 0$ if and only if $i+j-k =1 \lor k-\left|i-j\right| =1$. Note that it is actually now a natural option to include the domain via the parameter $a$ into the blow-up analysis. We shall see that the dynamic domain for the PDE nicely ``connects'' different regions in the desingularization procedure. In order to understand the behavior of~\eqref{eq:k0_general_forblowup} around the origin, we conduct a blow-up transformation with the following changes of coordinates, for $k \geq 1$:
\begin{equation}
u_k = \bar{r}^{\alpha_k} \bar{u}_k, \ v_k = \bar{r}^{\beta_k} \bar{v}_k, \ \varepsilon = \bar{r}^{\zeta} \bar \varepsilon, \
a = \bar{r}^{\delta} \bar a,
\end{equation}
where $\bar{r} \in [0, \rho]$, $\rho > 0$, $(\bar{u}_1, \bar{v}_1, \dots, \bar{u}_{k_0}, \bar{v}_{k_0}, \bar \varepsilon, \bar a)$ are coordinates in the ambient space of the manifold 
\benn
M = \mathbb{S}^{2k_0} \times [0, A]=\left\{\sum_{j=1}^{k_0}\left(\bar{u}_j^2+\bar{v}_j^2\right)+\bar \varepsilon^2=1\right\}\times [0, A],
\eenn 
where $\bar{a}\in [0,A]$ and $\alpha_k, \beta_k, \zeta, \delta \in \mathbb{R}$. The blow-up is a map from this cyclinder times a radial component into our original phase space
\benn
\Phi^{k_0}:M\times [0,\rho]\ra \R^{2k_0+2}.
\eenn
From the vector field $\cX_{k_0}$ induced by~\eqref{eq:k0_general_forblowup} one obtains via the pushforward condition $\Phi^{k_0}_*(\bar\cX_{k_0})=\cX_{k_0}$, a blown-up vector field $\bar\cX_{k_0}$. To determine the exponents defining the blow-up, one may use combinatorial theory based upon Newton polygons/polytopes~\cite{Bruno,Panazzolo2}, or directly exploit quasi-homogeneous scaling~\cite{Du93,KuehnBook}, whereas the latter method has turned out to be very successful for multiple time scale systems. Trying to obtain a quasi-homogeneous blow-up of polynomial order $N \in \mathbb{N}$, we get the following algebraic relations from equation~\eqref{eq:k0_general_forblowup}:
\begin{align}
u_1:  \, &N + \alpha_1 = 2 \alpha_1= 2 \beta_1 = \zeta + \delta,  \ \text{ and } N + \alpha_1 = \alpha_j + \beta_j, \, j \geq 2\,,\\
v_1: \,  &N + \beta_1 = \zeta + \delta\,, \\
u_k, k \geq 2: \, &N +\alpha_k = -\frac{3}{2} \delta + \alpha_k = \alpha_k + \alpha_1 = \beta_k + \beta_1, \\
 &N + \alpha_k = \alpha_j + \alpha_i = \beta_j + \beta_i,\, i,j \geq 2,\\
v_k, k \geq 2: \, &N + \beta_k \leq \zeta -\frac{3}{2} \delta + \beta_k \,.
\end{align}
We immediately observe that $-\frac{3}{2} \delta = N = \alpha_1 = \beta_1$ and $\alpha_k = \beta_k$ for all $k\geq 2$.
Furthermore, we observe that quasi-homogenity is not possible in $v_k$ unless $\zeta = 0$, which would imply $\delta = 2N = - 3 \delta$ and, hence, $ \delta = N = \alpha_k = \beta_k =0$, which is obviously not the desired situation. This is expected since the natural parabolic scaling of the heat operator part does not agree with the scaling enforced by the reaction terms. Therefore, we have to take $\zeta > 0$. Choosing the smallest integer coefficients, we obtain
\begin{equation}
N = \alpha_k = \beta_k = 3, k \geq 1, \ \delta = - 2, \ \zeta = 8.
\end{equation}
Note that $\delta$ being negative is a consequence of the fact that we assumed for the parameter of the domain $a=\cO(1)$ as $\varepsilon\ra 0$. Other scalings would certainly be favorable if we would deal with very large or very small domains; see Section~\ref{sec:conclusion}. The vector field $\bar\cX_{k_0}$ is smoothly conjugate to $\cX_{k_0}$ outside of the sphere $\mathbb{S}^{2k_0}$, yet new dynamics can be induced on the sphere, as we shall see below, by additional desingularization. Instead of directly working with the new coordinates, it is effective to parametrize the manifold $M\times [0,\rho]$ by different charts.

The fact that $\delta=-2$ means that in the following chart analysis we will have the relation $a = r_i^{-2} a_i$ for $i=1,2,3$. Hence, one has to be careful with the typical blow up investigation around $r_i=0$, since $r_i=0$, $a_i > 0$ does not correspond with original coordinates. 
However, an analysis is still feasible also in this case as all important perturbation arguments away from $r_i =0$ can still be lifted to the associated wedge for $r_i, a_i$, given by $a = r_i^{-2} a_i$.

\subsection{Analysis in first chart}

Building upon the ideas in~\cite{KruSzm4} and~\cite{EngelKuehn}, we analyze the dynamics through the origin by studying the system in three different charts. Firstly, we consider the \emph{entry chart} $K_1$, which is formally determined by setting $\bar u_1 = -1$ in the original blow-up. The precise coordinate transformation is given by
\begin{equation} \label{K1_scaling}
u_1 = -r_1^{3} , \ v_1 = r_1^{3}v_{1,1}, \
\varepsilon = r_1^8 \varepsilon_1, \ a = r_1^{-2} a_1, \    u_k = r_1^{3} u_{k,1}, \ v_k = r_1^{3} v_{k,1}, \ k\geq 2,
\end{equation}
where $r_1 \geq 0$. We can study the dynamics in this chart $K_1$ by applying the coordinate change~\eqref{K1_scaling}, identifying a common power of $r_1$, and then applying a suitable time change. This yields a desingularized vector field in the chart $K_1$. More precisely we have:  

\begin{proposition} 
\label{prop:K1}
Under the change of coordinates~\eqref{K1_scaling} system~\eqref{eq:k0_general}, for any $k_0 \in \mathbb{N}$, yields the desingularized blown-up system: 
\begin{subequations}
\label{eq:k0_K_1_general}
\begin{align}
\partial_{\tau} r_1&= - \frac{1}{3} r_1 F_1(r_1, v_{1,1}, \varepsilon_1, a_1, u_{2,1}, v_{2,1}, \dots)\,, \label{eq:r_1}\\
\partial_{\tau} \varepsilon_1&= \frac{8}{3} \varepsilon_1 F_1(r_1, v_{1,1}, \varepsilon_1, a_1, u_{2,1}, v_{2,1}, \dots)\,, \\
\partial_{\tau} v_{1,1}&= 2a_1\varepsilon_1  + v_{1,1} F_1(r_1, v_{1,1}, \varepsilon_1, a_1, u_{2,1}, v_{2,1}, \dots) + \mathcal{O}(r_1^{8}),\\
\partial_{\tau} u_{k,1}&= b_k a_1^{-3/2} u_{k,1} + 2 (-u_{k,1} - v_{k,1}v_{1,1}) 
+ \sqrt{2}\sum_{i,j=2}^{k_0} \alpha_{i,j}^k (u_{i,1} u_{j,1} - v_{i,1} v_{j,1}) \nonumber\\
&+ u_{k,1}F_1(r_1, v_{1,1}, \varepsilon_1, a_1, u_{2,1}, v_{2,1}, \dots) + \mathcal{O}(r_1^3), \ 2 \leq k \leq k_0, \\
\partial_{\tau} v_{k,1} &=  a_1^{-3/2} \left( r_1^{8}  b_k  v_{k,1} \varepsilon_1  +  a_1^{3/2}v_{k,1}F_1(r_1, v_{1,1}, \varepsilon_1, a_1, u_{2,1}, v_{2,1}, \dots) +  a_1^{3/2}\mathcal{O}(r_1^{8})\right), \ 2 \leq k \leq k_0,\\
\partial_{\tau} a_1&= - \frac{2}{3} a_1 F_1(r_1, v_{1,1}, \varepsilon_1, a_1, u_{2,1}, v_{2,1}, \dots)\,,
\end{align} 
\end{subequations}
where 
\begin{equation}
F_1(r_1, v_{1,1}, \varepsilon_1, a_1, u_{2,1}, v_{2,1}, \dots) = (1 - v_{1,1}^2) + 2a_1 \varepsilon_1 \mu + \sum_{j=2}^{k_0} u_{j,1}^2 - v_{j,1}^2 + \mathcal{O}(r_1^{3}).
\end{equation}
\end{proposition}
\begin{proof}
Firstly, we observe that
\begin{align*}
 3 r_1^2 \partial_t r_1 &= - \partial_t u_1 \\
 &= - \left( r_1^6 (1 - v_{1,1}^2) + r_1^6 2a_1 \varepsilon_1 \mu +  r_1^6 \sum_{j=2}^{k_0-1} u_{j,1}^2 - v_{j,1}^2 + \mathcal{O}\left(r_1^{9}\right)\right),
\end{align*}
and, hence, by a desingularization using a time change from $t$ to $\tau$ by eliminating the factor $r_1^3$, we obtain equation~\eqref{eq:r_1}. Similarly, we obtain 
$$\partial_{\tau} \varepsilon_1= - 8 \varepsilon_1 \frac{1}{r_1} \partial_{\tau} r_1 = \frac{8}{3} \varepsilon_1 F_1(r_1, v_{1,1}, \varepsilon_1, a_1, u_{2,1}, v_{2,1}, \dots).$$
Furthermore, we calculate
\begin{equation*}
3 r_1^2 v_{1,1} \partial_{\tau} r_1 + r_1^3 \partial_{\tau} v_{1,1} = \partial_{\tau} v_1 = 2a_1 r_1^3 \varepsilon_1 + \mathcal{O}\left(r_1^{11}\right) ,
\end{equation*}
which gives
$$ \partial_{\tau} v_{1,1} =  \frac{1}{r_1^3} \left(2a_1  r_1^3 \varepsilon_1 + v_{1,1} r_1^3 F_1(r_1, v_{1,1}, \varepsilon_1, a_1, u_{2,1}, v_{2,1}, \dots) + \mathcal{O}\left(r_1^{11}\right)\right).$$
The remaining computations are similar.
\end{proof}
For any $k_0 \geq 1$, note the existence of the three groups of equilibria
\begin{align} 
P_{\textnormal{a},1}^{-, k_0} &= \left\{ a_1 > 0\,:\,p_{\textnormal{a},1}^{-, k_0}(a_1) = \left(0,0,-1,0,0, \dots, a_1\right) \right\}, \label{equilibria_K1_ar_min} \\
P_{\textnormal{a},1}^{+, k_0} &=\left\{ a_1 > 0\,:\,p_{\textnormal{a},1}^{+, k_0}(a_1) = \left(0,0,1,0,0, \dots, a_1\right) \right\},\label{equilibria_K1_ar_plus}
\end{align}
and 
\begin{equation} \label{equilibria_K1_neutral}
Q_{1}^{\textnormal{in}, k_0} = \left\{a_1 > 0\,: \, q_{1}^{\textnormal{in}, k_0}(a_1) = (0,0,0,0,0 \dots,  a_1)\right\}.
\end{equation}
The equilibrium $q_{1}^{\textnormal{in}, k_0}(a_1)$ is a saddle with eigenvalue $1$ in the $v_{1,1}$-direction.

\begin{lemma} \label{lem:K1_equilibira}
The equilibria $p_{\textnormal{a},1}^{-, k_0}(a_1)$ and $p_{\textnormal{a},1}^{+, k_0}(a_1)$ are 
\begin{enumerate}
\item stable in the $v_{1,1}$-direction with eigenvalue $-2$,
\item stable in the $u_{k,1}$-directions with eigenvalues $-2+ b_k a_1^{-3/2}$,
\item  and neutral in all other directions.
\end{enumerate}
Furthermore, for $0 < r_1 < R_1$ sufficiently small, there are normally hyperbolic sets of equilibria $S_{\txta, 1}^{-, k_0}$ and $S_{\txta, 1}^{+, k_0}$ emanating from the lines $P_{\textnormal{a},1}^{-, k_0}$ and $P_{\textnormal{a},1}^{+, k_0}$.
\end{lemma}
\begin{proof}
The first three points follow from  straightforward calculations. The last statement can be directly deduced using the implicit function theorem.
\end{proof}
Note that $S_{\txta, 1}^{-,k_0} \cap \{ r_1^2 a = a_1\}$ coincides with the branch $S_{\txta}^{-,k_0}$ of the critical manifold for the Galerkin ODE~\eqref{eq:k0_general_forblowup} in original coordinates, and the same holds for $S_{\txta, 1}^{+, k_0} \cap \{ r_1^2 a = a_1\}$ and $S_{\txta}^{+, k_0}$. The relation $r_1^2 a = a_1$ means that $r_1 =0$ implies $a_1 =0$ if the direct correspondence to the original dynamics is preserved. Hence, one has to work in suitably restricted domains, e.g., for a given $a>0$ only consider in phase space all $a_1,r_1\geq 0$ such that $0< a_1/r_1^2<\I$, i.e., a suitable wedge-like region. Similar remarks apply to the other coordinates involving inverses such as $a_1^{-3/2}$. Henceforth, we make the following calculations with a general $a_1^* \approx 0$ within the allowed wedge-like domain, i.e., we are going to construct the center manifolds at $M_{\textnormal{a},1}^{-, k_0}(a_1^*)$ and $M_{\textnormal{a},1}^{+, k_0}(a_1^*)$ at the equilibria $p_{\textnormal{a},1}^{-, k_0}(a_1^*)$ and $p_{\textnormal{a},1}^{+, k_0}(a_1^*)$, allowing for an arbitrarily close approximation of the center manifolds at the desingularized origin.

\subsubsection{Center manifold approximation}
In the following, we will only focus on $p_{\textnormal{a},1}^{-, k_0}(a_1^*)$. The calculations for $p_{\textnormal{a},1}^{+, k_0}(a_1^*)$ will be analogous. Again, we will omit the higher order terms without loss of generality.

\subsubsection*{Example with three modes}
Firstly, we do the calculations for $k_0 =3$, in order to give more intuition on the relation between the higher modes and exemplifying the general calculations to follow.
For $k_0 =3$, equation~\eqref{eq:k0_K_1_general} (without terms of order $\mathcal{O}(r_1^3)$) reads
\begin{equation} \label{eq:k03_K1}
\begin{array}{lcl}
\partial_t v_{1,1}&=&  2a_1 \varepsilon_1 + v_{1,1} F_1,\\
\partial_t r_1&=& -\frac{1}{3} r_1 F_1, \\
\partial_t \varepsilon_1&=& \frac{8}{3} \varepsilon_1 F_1,\\
\partial_t a_1&=& -\frac{2}{3} a_1 F_1, \\
\partial_t u_{2,1}&=& b_2 a_1^{-3/2} u_{2,1} + 2 (-u_{2,1}  -  v_{2,1} v_{1,1})+  \sqrt{2}(u_{2,1} u_{3,1} -  v_{2,1} v_{3,1}) + u_{2,1} F_1,\\
\partial_t v_{2,1}&=& r_1^8 b_3 \varepsilon_1 a_1^{-3/2} v_{2,1} +  v_{2,1} F_1, \\
\partial_t u_{3,1}&=& b_3 a_1^{-3/2} u_{3,1} + 2 (-u_{3,1}  -  v_{3,1} v_{1,1})+  \frac{1}{\sqrt{2}}(u_{2,1}^2 -  v_{2,1}^2)+  u_{3,1} F_1,\\
\partial_t v_{3,1}&=& r_1^8 b_3 \varepsilon_1 a_1^{-3/2} v_{3,1} + v_{3,1} F_1,
\end{array}
\end{equation} 
where
$$ F_1 = (1- v_{1,1}^2) + 2 a_1 \varepsilon_1 \mu +  u_2^2 - v_2^2 + u_3^2 - v_3^2.$$
With a shift to $\tilde v_{1,1} = v_{1,1} +1$ and $ \tilde a_1 = a_1 - a_1^*$, we obtain
\begin{equation} \label{eq:k03_shifted_K1}
\begin{array}{lcl}
\partial_t \tilde v_{1,1}&=&  2(\tilde a_1 + a_1^*) \varepsilon_1 + (\tilde v_{1,1}-1) \tilde F_1,\\
\partial_t r_1&=& -\frac{1}{3} r_1 \tilde F_1, \\
\partial_t \varepsilon_1&=& \frac{8}{3} \varepsilon_1 \tilde F_1,\\
\partial_t \tilde a_1&=& -\frac{2}{3} \tilde a_1 \tilde F_1, \\
\partial_t u_{2,1}&=& b_2 (\tilde a_1+a_1^*)^{-3/2} u_{2,1} + 2 (-u_{2,1}  -  v_{2,1} (\tilde v_{1,1}-1))+   \sqrt{2}(u_{2,1} u_{3,1} -  v_{2,1} v_{3,1}) +  u_{2,1} \tilde F_1,\\
\partial_t v_{2,1}&=& r_1^8 b_3 \varepsilon_1 (\tilde a_1+a_1^*)^{-3/2} v_{2,1} +  v_{2,1} \tilde F_1, \\
\partial_t u_{3,1}&=& b_3 (\tilde a_1+a_1^*)^{-3/2} u_{3,1} + 2 (-u_{3,1}  -  v_{3,1} (\tilde v_{1,1}-1))+  \frac{1}{\sqrt{2}}(u_{2,1}^2 -  v_{2,1}^2)+  u_{3,1} \tilde F_1,\\
\partial_t v_{3,1}&=& r_1^8 b_3 \varepsilon_1 (\tilde a_1+a_1^*)^{-3/2} v_{3,1} +  v_{3,1} \tilde F_1,
\end{array}
\end{equation} 
where
$$ \tilde F_1 = (1- (\tilde v_{1,1}-1)^2) + 2 (\tilde a_1 + a_1^*) \varepsilon_1 \mu +  u_2^2 - v_2^2 + u_3^2 - v_3^2.$$
We now consider the equilibrium at the origin, with the three stable directions along $v_{1,1}$ with eigenvalue $\lambda_1 = -2$, $u_{2,1}$ with eigenvalue $\lambda_2 = -2 + b_2 (a_1^*)^{-3/2}$ and $u_{3,1}$ with eigenvalue $\lambda_3 = -2 + b_3 (a_1^*)^{-3/2}$, and otherwise center directions.  Transferring to standard form coordinates
\begin{align*}
\tilde v_{1,1} &= \frac{3}{2a_1^*}\left(y_1 + \frac{8}{9} (a_1^*)^3 \varepsilon_1 \left(\frac{3}{4a_1^*} - \frac{3 \mu}{4 a_1^*}\right)\right), \\
\tilde a_1 &= x_3 + \frac{2a_1^* \tilde v_{1,1}}{3} + \frac{8}{9} (a_1^*)^3 \varepsilon_1 \left(\frac{3}{4a_1^*} - \frac{3 \mu}{4 a_1^*}\right),\\
r_1 &= x_1,\ u_2 =y_2, \ u_3 = y_3 \
\varepsilon_1 = -3 \frac{x_2}{4 (a_1^*)^2},\
v_2 = x_4,\
v_3 = x_5,
\end{align*}
we consider
\begin{align*}
\hat h_1(x_1, x_2, x_3, x_4, x_5) &= b_{11} x_1^2 + b_{12} x_1 x_2 + b_{22} x_2^2 + b_{13} x_1 x_3 + b_{23} x_2 x_3 + b_{33} x_3^2 
+b_{44} x_4^2 + b_{14} x_1 x_4 \\
&+ b_{24} x_2 x_4 + b_{34} x_4 x_3 +b_{55} x_5^2 + b_{15} x_1 x_5 + b_{25} x_2 x_5 + b_{35} x_5 x_3 + b_{45} x_5 x_4 + \mathcal{O}(3), \\
\hat h_2(x_1, x_2, x_3, x_4) &= c_{11} x_1^2 + c_{12} x_1 x_2 + c_{22} x_2^2 + c_{13} x_1 x_3 + c_{23} x_2 x_3 + c_{33} x_3^2 
+c_{44} x_4^2 + c_{14} x_1 x_4 \\
&+ c_{24} x_2 x_4 + c_{34} x_4 x_3 +c_{55} x_5^2 + c_{15} x_1 x_5 + c_{25} x_2 x_5 + c_{35} x_5 x_3 + c_{45} x_5 x_4 + \mathcal{O}(3), \\
\hat h_3(x_1, x_2, x_3, x_4) &= d_{11} x_1^2 + d_{12} x_1 x_2 + d_{22} x_2^2 + d_{13} x_1 x_3 + d_{23} x_2 x_3 + d_{33} x_3^2 
+d_{44} x_4^2 + d_{14} x_1 x_4 \\
&+ d_{24} x_2 x_4 + d_{34} x_4 x_3 +d_{55} x_5^2 + d_{15} x_1 x_5 + d_{25} x_2 x_5 + d_{35} x_5 x_3 + d_{45} x_5 x_4 + \mathcal{O}(3).
\end{align*}
The invariance equation associated to~\eqref{eq:k03_shifted_K1} gives the coefficients
\begin{align} 
b_{22} &= \frac{3(\mu^2 -1)}{16 a_1^*}, \ b_{23} = \frac{\mu -1}{2 a_1^*}, \ b_{44} = \frac{a_1^*}{3}, \ b_{55} = \frac{a_1^*}{3}, \\
c_{24} &= \frac{3\sqrt{a_1^*}(\mu -1) }{2\left(b_2-2(a_1^*)^{3/2}\right)}, \ c_{45} = \frac{\sqrt{2}(a_1^*)^{3/2} }{b_2-2(a_1^*)^{3/2}} \\
d_{25}&= \frac{3\sqrt{a_1^*}(\mu -1) }{2\left(b_3-2(a_1^*)^{3/2}\right)}, \ d_{44} = \frac{(a_1^*)^{3/2} }{\sqrt{2}\left(b_2-2(a_1^*)^{3/2}\right)},
\end{align}
and all other coefficients being zero. Transforming back gives the center manifold approximation in form of the maps
\begin{align}
\tilde v_{1,1} &= \tilde h_1 (r_1, \varepsilon_1, \tilde a_1, v_{2,1}, v_{3,1}) \nonumber\\ 
&=\frac{1}{1 - \frac{2}{3} a_1^* \varepsilon_1 (\mu-1)} \left[ (1- \mu) (\tilde a_1 + a_1^*) \varepsilon_1 + (a_1^*)^2 \varepsilon_1^2 \frac{(\mu -1)(7 \mu -1)}{6} + v_{2,1}^2 +  v_{3,1}^2 \right] + \mathcal{O}(3), \label{eq:cmtilde_k03_1_K1} \\
 u_{2,1} &= \tilde h_2 (r_1, \varepsilon_1, \tilde a_1, v_{2,1}, v_{3,1}) \nonumber\\
 &= -\frac{2 (a_1^*)^{5/2}(\mu-1)}{b_2 - 2 (a_1^*)^{3/2}} v_{2,1} \varepsilon_1  + \frac{\sqrt{2}(a_1^*)^{3/2} }{b_2-2(a_1^*)^{3/2}}v_{2,1}v_{3,1} + \mathcal{O}(3), \label{eq:cmtilde_k03_2_K1} \\
 u_{3,1} &= \tilde h_2 (r_1, \varepsilon_1, \tilde a_1, v_{2,1}, v_{3,1}) \nonumber\\
 &= -\frac{2 (a_1^*)^{5/2}(\mu-1)}{b_3 - 2 (a_1^*)^{3/2}} v_{3,1} \varepsilon_1 + \frac{(a_1^*)^{3/2} }{\sqrt{2}\left(b_2-2(a_1^*)^{3/2}\right)} v_{2,1}^2 + \mathcal{O}(3). \label{eq:cmtilde_k03_3_K1}
\end{align}
One can check by computation that, indeed, $\tilde h= (\tilde h_1, \tilde h_2, \tilde h_3)$ satisfies the system of invariance equations associated with~\eqref{eq:k03_shifted_K1}. One can then transform back to $ v_{1,1}= \tilde v_{1,1} -1$ and $ a_1 = \tilde a_1 + a_1^*$.
\subsubsection*{General number of modes}
We now consider system~\eqref{eq:k0_K_1_general} with general $k_0 \in \mathbb{N}$, again shifting to $\tilde v_{1,1}^{k_0} = v_{1,1}^{k_0} +1$ and $ \tilde a_1 = a_1 - a_1^*$ to fix the equilibrium at the origin. We introduce
\begin{align}
\hat M_{\txta,1}^{-,k_0} = \bigg\{&\left( x_1^{k_0}, (r_1, \varepsilon_1, a_1, v_1^{k_0})\right) \in  \mathbb{R}^{k_0} \times \mathbb{R}^{k_0+2}\,:\, \nonumber\\
&\left( \tilde v_{1,1}^{k_0}, u_{2,1}^{k_0} \dots, u_{k_0,1}^{k_0} \right) =  h_1^{k_0}((r_1, \varepsilon_1, \tilde a_1, v_{2,1}^{k_0}, \dots, v_{k_0,1}^{k_0}) \bigg\}.\label{center_stable_K1_k0}
\end{align}
For $a_1^* \leq a_1 \leq A_1$, we introduce $\tilde A_1 = A_1 - a_1^*$.
For a neighbourhood $\mathcal{U}(0) \subset L^2([-\tilde A_1, \tilde A_1])$, , we might view the parametrization map $ h_1^{k_0}$ as a map
$$  h_1^{k_0} : [0, R_1] \times [0, \varepsilon_1^1] \times [0, \tilde A_1] \times \mathcal{U}(0) \to L^2([-\tilde A_1, \tilde A_1]), \ ( r_1, \varepsilon_1, a_1, v_1) \mapsto x_1^{k_0} $$
for some fixed $R_1, \varepsilon_1^1, \tilde A_1 > 0$, via identifying
\begin{align*}
 v_1 &=  \sum_{k=2}^{\infty} v_{k,1} e_k, \\
 x_1^{k_0}  &= \tilde v_{1,1}^{k_0} e_1 + \sum_{k=2}^{k_0} u_{k,1}^{k_0} e_k,
\end{align*}
and denoting
$ v_1^{k_0} = \left(0, v_{2,1}, \dots,v_{k_0,1} \right)$.
Hence, we can consider $\hat M_{\txta,1}^{-,k_0}$ also as a Banach manifold
\begin{align} 
 M_{\txta,1}^{-,k_0} = \bigg\{ \left( x_1, (r_1, \varepsilon_1, \tilde a_1, v_1)\right) \in  & L^2([-\tilde A_1, \tilde A_1]) \times \left( [0, R_1] \times [0, \varepsilon_1^1] \times [0, \tilde A_1] \times \mathcal{U}(0) \right) \,:\, \nonumber\\
& x_1 =  h_1^{k_0}(r_1, \varepsilon_1, \tilde a_1, v_1) \bigg\}. \label{center_stable_K1_k0_L2}
\end{align}

Similarly to Lemma~\ref{lem:solution_of_invariance_equation}, we can show the following:
\begin{lemma} \label{lem:solution_of_invariance_equation_K1}
The invariance equation associated to system~\eqref{eq:k0_K_1_general} --- under change of coordinates $\tilde v_{1,1}^{k_0} = v_{1,1}^{k_0} +1$ and $ \tilde a_1 = a_1 - a_1^*$ to shift the equilibrium to the origin --- has the solution
\begin{align}
\tilde v_{1,1}^{k_0} &= h_{1,1}^{k_0}(r_1,\varepsilon_1, \tilde a_1, v_1^{k_0}) \nonumber
\\ &= \frac{1}{1 - \frac{2}{3} a_1^* \varepsilon_1 (\mu-1)} \left[ (1- \mu) (\tilde a_1 + a_1^*) \varepsilon_1 + (a_1^*)^2 \varepsilon_1^2 \frac{(\mu -1)(7 \mu -1)}{6}  + \sum_{k=2}^{k_0} v_{k,1}^2 \right] + \mathcal{O}(3), \label{eq:cmtilde_k0gen_1_K1} \\
u_{k,1}^{k_0} &= h_{k,1}^{k_0}(r_1,\varepsilon_1, \tilde a_1, v_1^{k_0})  \nonumber \\
&= -\frac{2 (a_1^*)^{5/2}(\mu-1)}{b_k - 2 (a_1^*)^{3/2}} v_{k,1} \varepsilon_1 + C_k \frac{(a_1^*)^{3/2}}{(b_k - 2 (a_1^*)^{3/2})} \sum_{i,j=2}^{k_0} \beta_{i,j}^k v_i v_j +  \mathcal{O}(3), \label{eq:cmtilde_k0gen_k_K1}
\end{align}
where $ \left| C_k \right| < C$ for a constant $ C \in \mathbb{R}$ and $ \left| \beta_{i,j}^k \right| \in [0,1]$ with $\beta_{i,j}^k \neq 0$ if and only if $i+j-k =1 \lor k-\left|i-j\right| =1$. 
\end{lemma}
\begin{proof}
Follows from direct calculations similarly to Lemma~\ref{lem:solution_of_invariance_equation}.
\end{proof}

\begin{remark}
We can directly see the connection to the center manifold approximations in \cite{KruSzm4, EngelKuehn}, where only the parametrization of the $v_{1,1}$-direction is relevant. Note that on the invariant set $\{r_1 = u_{j,1} = v_{j,1} = 0, j\geq 2\}$, system~\eqref{eq:k0_K_1_general} becomes
\begin{align*}
\partial_{\tau} v_{1,1}&= 2 a_1 \varepsilon_1 +  v_{1,1}  F_1(0,v_{1,1}, \varepsilon_1,a_1, \dots),\\
\partial_{\tau} \varepsilon_1 &= \frac{8}{3} \varepsilon_1 F_1(0,v_{1,1}, \varepsilon_1,a_1, \dots),\\
\partial_{\tau} a_1 &= -\frac{2}{3} \varepsilon_1 F_1(0,v_{1,1}, \varepsilon_1,a_1, \dots).
\end{align*} 
Note that one can summarize $\rho_1 = 2 a_1 \varepsilon_1$ such that the above equation becomes
\begin{align*}
\partial_{\tau} v_{1,1}&= 2 \rho_1 +  v_{1,1}  F_1(0,v_{1,1}, \rho_1, \dots),\\
\partial_{\tau} \rho_1 &= 2 \rho_1 F_1(0,v_{1,1}, \rho_1, \dots).
\end{align*} 
We define
\begin{align}
h^-(\rho_1) &= -1 + (1-\mu) \rho_1 + \mathcal O\left(\rho_1^2\right),\\
h^+(\rho_1) &= 1 + (1+\mu) \rho_1 + \mathcal O\left(\rho_1^2\right).
\end{align}
Observe that
$$   F_1(0,h^-(\rho_1), \rho_1, \dots) = 2 \rho_1 + \mathcal O\left(\rho_1^2\right), \quad  F_1(0,h^+(\rho_1), \rho_1, \dots) = - 2 \rho_1 + \mathcal O\left(\rho_1^2\right).
$$
and, hence, we observe that the invariance equation is satisfied:
\begin{align*}
\mathcal O\left(\rho_1^2\right) = 2 \rho_1 + h^-(\rho_1) \left(2 \rho_1 + \mathcal O\left(\rho_1^2\right) \right) = \partial_{\tau}(h^-(\rho_1))  &= \frac{\partial h^-(\rho_1)}{\partial \rho_1} \partial_{\tau} \rho_1  =  \mathcal O\left(\rho_1^2\right),\\
\mathcal O\left(\rho_1^2\right) = 2 \rho_1 + h+(\rho_1) \left(-2 \rho_1 + \mathcal O\left(\rho_1^2\right) \right) = \partial_{\tau}(h^+(\rho_1))  &= \frac{\partial h^+(\rho_1)}{\partial \rho_1} \partial_{\tau} \rho_1  =  \mathcal O\left(\rho_1^2\right).
\end{align*}
Hence, for the reduced system on $\{r_1 = u_{j,1} = v_{j,1} = 0, j\geq 2\}$ we obtain the center manifold approximation
\begin{equation} \label{eq:Hpm}
H^{\pm} (\varepsilon_1, a_1, 0, 0, \dots ) = (1 \pm \mu) a_1 \varepsilon_1 \pm 1  
+ \mathcal{O}\left( a_1^2\varepsilon_1^2 \right).
\end{equation}
In equation~\eqref{eq:cmtilde_k0gen_1_K1}, setting $a_1^*=0$  and shifting $\tilde v_{1,1}^{k_0} = v_{1,1}^{k_0} +1$, yields exactly the same formula as for $H^-$. For $H^+$, the situation is analogous.
\end{remark}
Furthermore, we can deduce the following proposition:
\begin{proposition} \label{prop:convergence_K_1}
For any $v_{1} \in \mathcal{U}(0) \subset L^2([-\tilde A_1,\tilde A_1])$,  $0 \leq \tilde a_1 \leq \tilde A_1$ and $\varepsilon_1 \in [0, \varepsilon_1^1]$, $r_1 \in [0, R_1]$ small, we have that $\tilde v_{1,1}^{k_0}$ is uniformly bounded for all $k_0$ up to $\mathcal{O}(3)$ terms, and
\begin{align}
\left|u_{k,1}^{k_0}\right|^2 &= \left|h_k^{k_0}(r_1,\varepsilon_1, a_1, v_{1}^{k_0}) \right|^2   \leq \frac{C_1}{\left|b_k - 2 (a_1^*)^{3/2}\right|} (\|v_1\|_2^2 + \varepsilon_1) + \mathcal{O}(3) \label{eq:bound_uk_K1}
\end{align}
for all $k_0 \in \mathbb{N}$ and a constant $C_1 >0$.
\end{proposition}
\begin{proof}
Inequality~\eqref{eq:bound_uk_K1} follows immediately from Lemma~\ref{lem:solution_of_invariance_equation_K1} with similar arguments as in  Remark~\ref{prop:convergence_beforeblowup}.\end{proof}
We obtain the following corollary, where the manifolds $M_{\txta,1}^{-, n}(a_1^*)$ can be naturally transformed back via $ - 1 + \tilde v_{1,1}^{k_0} = v_{1,1}^{k_0}$ and $ \tilde a_1 + a_1^* = a_1$ to the corresponding objects for $ L^2([-A_1,A_1])$ and $a_1^* \leq a_1 \leq A_1$:
\begin{corollary} \label{cor:Maconvergence}
The manifolds $M_{\txta,1}^{-, n}(a_1^*)$, as given in~\eqref{center_stable_K1_k0},
 exist for all $a_1^* \geq 0$, in particular for $a_1^*=0$.
\end{corollary}
\begin{proof}
The first statement can be deduced directly from formulas~\eqref{eq:cmtilde_k0gen_1_K1} and~\eqref{eq:cmtilde_k0gen_k_K1} in Lemma~\ref{lem:solution_of_invariance_equation_K1}, which do not exhibit a singularity at $a_1^*=0$. 
\end{proof}
Note that this corollary allows us to understand the center manifolds $M_{\txta,1}^{-,n}(0)$ truly from the origin expressed by $r_1=0$, still respecting the relation $r_1^2 a = a_1$ enforced at some fixed $a>0$. In particular, the center manifolds exists within a suitable open set in $\{a_1\geq 0,r_1\geq 0\}$.

Similarly to Remark~\ref{rem:convergence_beforeblowup}, if we disregard the $\mathcal{O}(3)$ terms, we obtain that
 $$h_1^{n}(r_1, \varepsilon_1, \tilde a_1, v_{1}) = x_{1}^{n} \to x_{1}^*=:h_1(r_1, \varepsilon_1, \tilde a_1, v_{1} ) $$
as $n \to \infty$ in $L^2([-\tilde A_1,\tilde A_1])$, uniformly in $(r_1, \varepsilon_1, \tilde a_1, v_{1})$.
Analagously to~\eqref{center_stable_before_blowup}, one could then define the infinite-dimensional Banach manifold
\begin{align} 
 M_{\txta,1}^{-}(a_1^*) := \bigg\{&\left( x_1, (r_1, \varepsilon_1, \tilde a_1, v_1)\right) \in  L^2([-\tilde A_1, \tilde A_1]) \times \left(\mathbb{R}^3 \times  L^2([-\tilde A_1,\tilde A_1]) \right)\,:\, \nonumber\\
 & x_1 = h_1(r_1, \varepsilon_1, \tilde a_1, v_{1} ) \bigg\}.\label{center_stable_K1}
\end{align}
The manifolds $M_{\txta,1}^{-, n}(a_1^*)$, as given in~\eqref{center_stable_K1_k0},
converge to $M_{\txta,1}^{-}(a_1^*)$~\eqref{center_stable_K1} as $n \to \infty$ in Hausdorff distance with respect to the $L^2$-norm, if it exists. This can be seen from Proposition~\ref{prop:convergence_K_1} in an analogous manner to Remark~\ref{rem:Mcconvergence}.
However, we recall the caveat that such a limiting procedure is most likely not admissible.

\subsubsection*{Coincidence of manifolds}
In the following, we make a crucial observation for the use of blow-up in our Galerkin scheme, concerning the invariant manifolds 
$M_c^n$ before blow-up, the global blow-up manifolds 
$\bar M_a^{-,n}$
and its $K_1$ versions 
$M_{\txta,1}^{-,n}(0)$, and similarly for 
$ \bar M_a^{+,n} $ 
and its $K_1$ versions 
$ M_{\txta,1}^{+,n}(0)$. 
We denote, for each $n \in \mathbb{N}$, the blow-up map by $\Phi^n$, which is a diffeomorphism away from the origin, and the chart-$K_1$ map by $\kappa_1^n$.
\begin{proposition} \label{prop:manifold_convergence}
For all $ n \in \mathbb{N}$ and $c < 0$ sufficiently small, we have that  $M_c^n = \Phi^n(\bar M_a^{-,n})$ in a small neigbourhood of $\{u=v=c,\, \varepsilon =0\}$ and $M_{\txta,1}^{-,n}(0) = \kappa_1 (\bar M_a^{-,n})$, and $\tilde M_c^n = \Phi^n(\bar M_a^{+,n})$ in a small neigbourhood of $\{u=-v=c,\, \varepsilon =0\}$ and $M_{\txta,1}^{+,n}(0) = \kappa_1 (\bar M_a^{+,n})$.
\end{proposition}
\begin{proof}
The statement follows directly from the fact that we are dealing with smooth coordinate transformations of vector fields governing ODEs. 
\end{proof}
In the following, we will simply write 
$$ M_{\txta,1}^{-,n} := M_{\txta,1}^{-,n}(0), 
\ \text{ and } M_{\txta,1}^{+,n} := M_{\txta,1}^{+,n}(0).
$$
\subsubsection{Corresponding PDE}
Additionally, we can show the following proposition, concerning the relation to a PDE in chart $K_1$ with a free boundary:
\begin{proposition}
\label{prop:PDEcorrespond}
Keeping the dynamics in $\varepsilon_1, a_1, r_1$ the same, system~\eqref{eq:k0_K_1_general} is --- up to time changes --- the Galerkin discretization of the PDE system 
\be
\label{eq:mainPDE_K1_fullsystem}
\begin{array}{lcl}
\partial_{\tau} u_1&=&\sqrt{2 a_1}\left(\partial_x^2 u_1 + u_1^2-v_1^2 \right) + F(\varepsilon_1, r_1, a_1, u_1, v_1)u_1 + \mathcal{O}\left(r_1^3\right),\quad u_1=u_1(x,t),~v_1=v_1(x,t),\\
\partial_{\tau} v_1&=& \sqrt{2 a_1}\left(\varepsilon_1 r_1^8 \partial_x^2 v_1 + \varepsilon_1\right) + F(\varepsilon_1, r_1, a_1, u_1, v_1)v_1 + \mathcal{O}\left( r_1^{8}\right), \quad ~(x,t)\in[-a_1,a_1]\times[0,T], \\
\partial_{\tau} r_1&=& - \frac{1}{3} r_1 F(\varepsilon_1, r_1, a_1, u_1, v_1)\,,\\
\partial_{\tau} \varepsilon_1&=& \frac{8}{3} \varepsilon_1 F(\varepsilon_1, r_1, a_1, u_1, v_1)\,, \\
\partial_{\tau} a_1&=& - \frac{2}{3} a_1 F(\varepsilon_1, r_1, a_1, u_1, v_1),
\end{array} 
\ee
for solutions $(u_1,v_1)$ such that $u_{1,1} \equiv -1$, where the interval $[-a_1, a_1]$ is moving in time, and 
$F$ is given by
$$ F(\varepsilon_1, r_1, a_1, u_1, v_1) =  2a_1 \varepsilon_1 \mu + \left(\|u_1\|^2 - \|v_1\|^2 + \mathcal{O}\left( r_1^3\right) \right).$$
\end{proposition}
\begin{proof}
Firstly, observe that for all $k_0 > 0$,
\begin{align*}
&\left\langle \left( \sum_{j=1}^{k_0} u_{j,1} e_j \right)^2 - \left( \sum_{j=1}^{k_0} v_{j,1} e_j \right)^2, e_k\right\rangle \\ &= \sum_{j=1}^{k_0} u_{j,1}^2 \langle e_j^2, e_k\rangle + \sum_{j,i=1, j\neq i}^{k_0} u_{j,1} u_{i,1} \langle e_j e_i, e_k \rangle - \left( \sum_{j=1}^{k_0} v_{j,1}^2 \langle e_j^2, e_k\rangle + \sum_{j,i=1, j\neq i}^{k_0} v_{j,1} v_{i,1} \langle e_j e_i, e_k \rangle \right),
\end{align*}
and the same for $k_0 \to \infty$.
Secondly, we know
\begin{align*}
\lambda_1 &= 0, \ \langle 1, e_1 \rangle = \sqrt{2a_1}, \\
\langle e_k e_l, e_1 \rangle &= \langle e_1 e_l, e_k \rangle = \frac{1}{\sqrt{2a_1}}  \delta_{k,l}, \quad k,l \geq 1.
\end{align*}
Hence, we obtain, in addition to the equations for $\varepsilon_1, a_1, r_1$ as already given by~\eqref{eq:k0_K_1_general},
\begin{align*}
\partial_{\tau} v_{1,1}&= 2a_1\varepsilon_1  + v_{1,1} \left((1 - v_{1,1}^2) + 2a_1 \varepsilon_1 \mu +\sum_{j=2}^{\infty} u_{j,1}^2 - v_{j,1}^2 + \mathcal{O}\left( r_1^3\right) \right) + \mathcal{O}\left( r_1^{8}\right),\\
\partial_{\tau} u_{k,1}&= b_k a_1^{-3/2} u_{k,1} + 2 (-u_{k,1} - v_{k,1}v_{1,1}) 
+ \sqrt{2}\sum_{i,j=2}^{\infty} \alpha_{i,j}^k (u_{i,1} u_{j,1} - v_{i,1} v_{j,1}) \nonumber\\
&+  u_{k,1} \left((1 - v_{1,1}^2) + 2a_1 \varepsilon_1 \mu +\sum_{j=2}^{\infty} u_{j,1}^2 - v_{j,1}^2 + \mathcal{O}\left( r_1^3\right) \right)+ \mathcal{O}\left( r_1^3\right) , \ k \geq 2, \\
\partial_{\tau} v_{k,1} &=   r_1^{8}  b_k a_1^{-3/2} v_{k,1} \varepsilon_1  +  v_{k,1}\left((1 - v_{1,1}^2) + 2a_1 \varepsilon_1 \mu +\sum_{j=2}^{\infty} u_{j,1}^2 - v_{j,1}^2 + \mathcal{O}\left( r_1^3\right)\right) + \mathcal{O}\left( r_1^{8}\right), \ k \geq 2.
\end{align*}
Now these equations can be truncated at any $k_0 \geq 1$ as before, which shows the claim.
\end{proof}

\begin{figure}[ht]
        \centering
        \begin{subfigure}{0.5\textwidth}
        \centering
  		\begin{overpic}[width=1\linewidth]{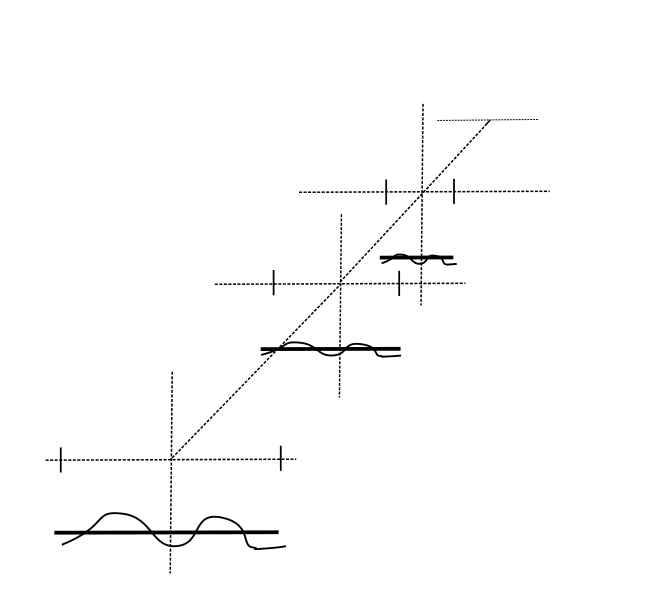}   
\put(46,12){\tiny $v_{1,1}\approx-1$}	
	\put(17,35){\tiny $x_1(x)$}
	\put(43,58){\tiny $x_1(x)$}
	\put(55,76){\tiny $x_1(x)$}
  		\put(8,20){\tiny $-a_1$}
  		\put(41,20){\tiny $a_1$}
  		\put(55,61){\tiny $-a_1$}
  		\put(67,61){\tiny $a_1$}
  		\put(38,47){\tiny $-a_1$}
  		\put(59,47){\tiny $a_1$}
  		\put(84,75){\small $r_1 =0$}
                \put(52,22){\small $r_1 = R_1$}
        \end{overpic}
        \caption{Extracts around $M_{\txta,1}^{-, k_0}$}
        \label{Mminus}
		\end{subfigure}%
        \begin{subfigure}{0.5\textwidth}
        \centering
  		\begin{overpic}[width=1\linewidth]{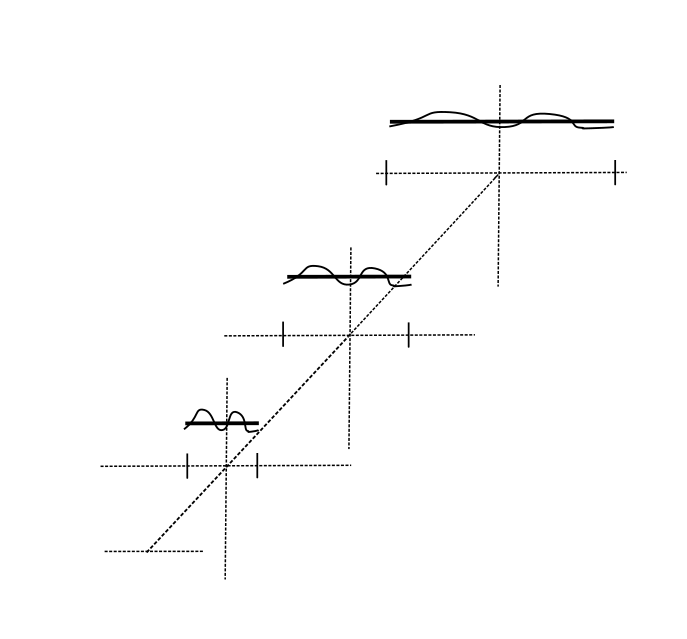}   
\put(12,31){\tiny $v_{1,1} \approx 1$}	
	\put(25,39){\tiny $x_1(x)$}
	\put(43,58){\tiny $x_1(x)$}
	\put(65,83){\tiny $x_1(x)$}
  		\put(22,21){\tiny $-a_1$}
  		\put(37,21){\tiny $a_1$}
  		\put(54,64){\tiny $-a_1$}
  		\put(88,64){\tiny $a_1$}
  		\put(38,41){\tiny $-a_1$}
  		\put(59,41){\tiny $a_1$}
  		\put(37,68){\small $r_1 =R_1$}
                \put(35,12){\small $r_1 = 0$}
        \end{overpic}
        \caption{Extracts around $M_{\txta,1}^{+, k_0}$}
        \label{Mplus}
		\end{subfigure}
		\caption{Behavior of solutions $x_1= v_{1,1} e_1 + \sum_{k=2}^{k_0} u_{k,1}^{k_0} e_k$ for the Galerkin problem in chart $K_1$ along $M_{\txta,1}^{-, k_0}$ (a), with decreasing $a_1$ and $r_1$, and $M_{\txta,1}^{+, k_0}$ (b), with increasing $a_1$ and $r_1$. Note that replacing $v_{1,1}$ by $u_{1,1} =-1$ means replacing $x_1$ by $u_1$, approixmating the solution of the corresponding PDE in~\eqref{eq:mainPDE_K1_fullsystem}. The behavior around $r_1=0$, i.e.~in between the two illustrated regions, is captured by the analysis in chart $K_2$.}
        \label{fig:K1}
\end{figure}

The importance of Proposition~\ref{prop:PDEcorrespond} is twofold. On the one hand, there is a free boundary showing that adjusting the spatial domain dynamically is helpful in a PDE blow-up. On the other hand, we observe that the PDE is quasi-linear, yet still locally well-posed, as the diffusion coefficients are dynamic variables. Both aspects are completely new and cannot appear for ODEs; see also Section~\ref{sec:conclusion}. The correspondence between Proposition~\ref{prop:convergence_K_1} and Corollary~\ref{cor:Maconvergence} on the one side and Proposition~\ref{prop:PDEcorrespond} on the other is illustrated in Figure~\ref{fig:K1}.

\subsection{Analysis in second chart} 
\label{sec:chart2}
We turn to analyzing the system in the rescaling chart $K_2$, determined formally by setting $\bar \varepsilon = 1$, where the dynamics around the origin can be understood in detail. The system of Galerkin ODEs in $K_2$ is given by the following result:

\begin{proposition}
Consider the rescaling
\begin{equation} \label{K2_scaling}
u_1 = r_2^3 u_{1,2}, \ v_1 = r_2^3 v_{1,2}, \ u_k = r_2^{3} u_{k,2}, \ v_k = r_2^{3} v_{k,2}, \
\varepsilon = r_2^8, \ a = r_2^{-2} a_2,
\end{equation}
where $r_2 \geq 0$.
Under this change of coordinates and desingularization, system~\eqref{eq:k0_general} becomes, for any $k_0 \geq 1$,
\begin{subequations}
\label{eq:k0_K_2_general}
\begin{align}
\partial_{\tau} u_{1,2}&= u_{1,2}^2 - v_{1,2}^2 + 2a_2 \mu +  \sum_{j=2}^{k_0} \left(u_{j,2}^2 - v_{j,2}^2 \right) + \mathcal{O} \left(r_2^3 \right),\\
\partial_{\tau} v_{1,2}&= 2a_2 + \mathcal{O} \left(r_2^{8} \right), \\
\partial_{\tau} u_{k,2}&= b_k a_2^{-3/2} u_{k,2} + 2 (u_{k,2} u_{1,2}  -  v_{k,2} v_{1,2}) \nonumber\\
&+ \sqrt{2}\sum_{i,j=2}^{k_0} \alpha_{i,j}^k (u_{i,2} u_{j,2} - v_{i,2} v_{j,2}) + \mathcal{O} \left(r_2^3 \right), \ 2 \leq k\leq k_0,\\
\partial_{\tau} v_{k,2}&= a_2^{-3/2} \left( r_2^8 b_k v_{k,2} + a_2^{3/2}\mathcal{O} (r_2^{8}) \right), \ 2 \leq k \leq k_0, \label{eq:k0_K_2_general_vk} \\
\partial_{\tau} a_2 &= 0, \\
\partial_{\tau} r_2 &= 0.
\end{align} 
\end{subequations}
\end{proposition}
\begin{proof}
Clearly, we have $\partial_t r_2 = 0$ and also $\partial_{\tau} r_2 = 0$ for $\tau = r_2^3 t$. Hence, we observe that
\begin{equation*}
r_2^3 \partial_t u_{1,2} = \partial_t u_1 =  r_2^6 \left( u_{1,2}^2 -  v_{1,2}^2\right) + r_2^8 r_2^{-2} 2 a_2 \mu + \sum_{j=2}^{k_0} r_2^{6} \left( u_{j,2}^2 -  v_{j,2}^2 \right) + \mathcal O \left(r_2^9\right),
\end{equation*}
which gives
$$ \partial_t u_{1,2} = r_2^3 \left( u_{1,2}^2 -  v_{1,2}^2\right) + r_2^3 2 a_2 \mu + r_2^3 \sum_{j=2}^{k_0 -1} \left( u_{j,2}^2 -  v_{j,2}^2 \right)+ \mathcal O \left(r_2^6\right) .$$
Hence, by the time rescaling $\tau = r_2^3 t$, we obtain the claimed equation. Additionally, we compute
$$ r_2^3 \partial_t v_{1,2} = \partial_t v_1 = r_2^6 2 a_2 + \mathcal O \left(r_2^{14}\right),$$
such that $\partial_{\tau} v_{1,2} =  2a_2  + \mathcal O \left(r_2^{8}\right)$ follows. The remaining equations can be deduced similarly.
\end{proof}

We can also compute the transition maps between $K_1$ and $K_2$:
\begin{lemma}
The changes of coordinates between the charts $K_1$ and $K_2$ for $k_0 \in \mathbb{N}$ are given by the maps
$\kappa_{12}^{k_0}: K_1 \to K_2$ with
\begin{equation} \label{kappa12}
u_{1,2} = - \varepsilon_{1}^{-3/8}, \ v_{1,2} = \varepsilon_{1}^{-3/8} v_{1,1}, \ u_{k,2} = \varepsilon_1^{-3/8} u_{k,1}, \ v_{k,2} = \varepsilon_1^{-3/8} v_{k,1}, \ r_2 = \varepsilon_1^{1/8} r_1, \ a_2 = \varepsilon_1^{1/4} a_1,
\end{equation}
and
$\kappa_{21}^{k_0}: K_2 \to K_1$ with
\begin{equation} \label{kappa21}
\varepsilon_{1} = u_{1,2}^{-8/3}, \ v_{1,1} = -u_{1,2}^{-1} v_{1,2}, \ u_{k,1} = -u_{1,2}^{-1} u_{k,2}, \ v_{k,1} = -u_{1,2}^{-1} v_{k,2}, \ r_1 = - u_{1,2}^{1/3} r_2, \ a_1 = u_{1,2}^{2/3} a_2.
\end{equation}
\end{lemma}
\begin{proof}
Follows from straight-forward calculations.
\end{proof}
Note that for $r_2 = 0$, equation~\eqref{eq:k0_K_2_general} becomes
\begin{align*}
\partial_{\tau} u_{1,2}&= \left(u_{1,2}^2 -  v_{1,2}^2\right) + 2 a_2 \mu + \sum_{j=2}^{k_0} \left(u_{j,2}^2 - v_{j,2}^2 \right),\\
\partial_{\tau} v_{1,2} &=  2a_2, \\
\partial_{\tau} u_{k,2}&= - \sqrt{2} \pi^2 (k-1)^2 a_2^{-3/2} u_{k,2} + 2 (u_{k,2}u_{1,2}-v_{k,2}v_{1,2}) + \sqrt{2}\sum_{i,j=2}^{k_0} \alpha_{i,j}^k (u_{i,2} u_{j,2} - v_{i,2} v_{j,2}), \ k \geq 2, \\
\partial_{\tau} v_{k,2} &=   0, \ \partial_{\tau} r_2 = 0, \ \partial_{\tau} a_2 = 0.
\end{align*}
Taking $v_{k,2}$ zero for $k \geq 2$, we observe that, at $u_{k,2} =0$, the linearization has the factor
$$  \left( 2 u_{1,2}  - \sqrt{2}(k-1)^2 \pi^2 a_2^{-3/2} \right) u_{k,2}.$$
Hence, $u_{k,2} =0$, for all $k\geq0$, is an exponentially stable equilibrium as long as 
\begin{equation} \label{eq:u12_bound}
u_{1,2} < \frac{\pi^2}{\sqrt{2}} \frac{1}{a_2^{3/2}} < (k-1)^2\frac{\pi^2}{\sqrt{2}} \frac{1}{a_2^{3/2}} \quad \forall k \geq 2.
\end{equation}
Hence, for $u_{1,2}, v_{1,2}$ the analysis of the transcritical singularity, depending on the sign of $\mu$, can be built upon similar ideas for the ODE case~\cite{EngelKuehn, KruSzm4}, which covers only the first mode in our context. We remark that~\cite{KruSzm4} relies on a graphical phase plane argument in the scaling chart $K_2$, so for analytic estimates of the dynamics in the scaling chart $K_2$ we refer to \cite{EngelKuehn}, which treats the Euler discretization of the normal form ODE but, by taking the discretization parameter $h \to 0$, provides a detailed analysis of the dynamics also for the ODE case, i.e., for the first mode in our setting. Building upon these previous results, we are going to prove the main dynamical feature for equation~\eqref{eq:k0_K_2_general} with 
arbitrarily large $k_0 \in \mathbb N$
in Proposition~\ref{propK2}. 
Note that in Appendix~\ref{ap:alternative_proof}, we sketch an alternative direct PDE argument to the following proof of Proposition~\ref{propK2}.
 
For suitably small $\delta, \beta, \nu > 0$, the following sets are understood as subsets of $\R^{k_0} \times \R^{k_0} \times [0, A_2] \times [0, R_2]$,
where $A_2 = \delta^{-1/4} A_1, R_2 = \delta^{1/8} R_1$ are inherited from the first chart;
we set
\begin{align*}
 \Sigma_2^{\textnormal{in}, k_0} = \bigg\{ u_{1,2} = - \delta^{-3/8},& v_{1,2} \in (- \delta^{-3/8} - \beta, - \delta^{-3/8} + \beta), a_2 =\delta^{1/4} \nu, \\
 & u_{k,2}^2 \leq C_{u, k_0} \frac{\delta^{3/4}}{\left|b_k\right|} , v_{k,2}^2 \leq C_{v, k_0} \frac{\delta^{3/4}}{\left|b_k\right|}, \, 2 \leq k \leq k_0  \bigg\},
\end{align*}
$$ \Sigma_{2,\txta }^{\textnormal{out}, k_0} = \left\{ u_{1,2} = - \delta^{-3/8}, v_{1,2} \in ( \delta^{-3/8} - \beta,  \delta^{-3/8} + \beta), \, u_{k,2}^2 \leq C_{u, k_0} \frac{\delta^{3/4}}{\left|b_k\right|} , v_{k,2}^2 \leq C_{v, k_0} \frac{\delta^{3/4}}{\left|b_k\right|}, \, 2 \leq k \leq k_0 \right\}, $$ 
and 
$$\Sigma_{2,\txte }^{\textnormal{out}, k_0} = \left\{ u_{1,2} =  \delta^{-3/8}, v_{1,2} \leq  \Omega(\mu)\delta^{-1/8}, \, u_{k,2}^2 \leq C_{u, k_0} \frac{\delta^{3/4}}{\left|b_k\right|} , v_{k,2}^2 \leq C_{v, k_0} \frac{\delta^{3/4}}{\left|b_k\right|}, \, 2 \leq k \leq k_0  \right\},$$
where $\Omega(\mu) > 0$.
Similarly to Proposition~\ref{prop:manifold_convergence}, 
we have
\begin{align*}
M_{\txta,2}^{-,k_0} &= \kappa_{12}^{k_0}(M_{\txta,1}^{-, k_0}),  
\end{align*} 
Note with inequality~\eqref{eq:bound_uk_K1} that 
for appropriate $C_{u, k_0}, C_{v, k_0}, \beta > 0$ we have that, for sufficiently small $v_{k,2}$, $$\kappa_{12}^{k_0}(M_{\txta,1}^{-, k_0}) \cap \{u_{1,2} = - \delta^{-3/8} \} = M_{\txta,2}^{-, k_0} \cap \{u_{1,2} = - \delta^{-3/8} \} \subset \Sigma_2^{\textnormal{in}, k_0}.$$
In accordance with this, we choose
$$\Sigma_1^{\textnormal{out}, k_0} = \left\{ \varepsilon_1 = \delta, a_1 = \nu,  u_{k,1}^2 \leq C_{u, k_0} \frac{\varepsilon_1^{3/2}}{\left|b_k\right|} , v_{k,1}^2 \leq C_{v, k_0} \frac{\varepsilon_1^{3/2}}{\left|b_k\right|}, 2 \leq k \leq k_0 \right\},$$
and 
$$ \Sigma_1^{\textnormal{in}, k_0} = \{ r_1 = \rho \}.$$
Note that we will later, in the proof of Theorem~\ref{thm:mainresult}, introduce a small subset $W_1^{k_0} \subset \Sigma_1^{\textnormal{in}, k_0}$ which forms a neighbourhood of $M_{\txta,1}^{-, k_0} \cap \Sigma_1^{\textnormal{in}, k_0}$ for sufficiently small $v_{k,1}$ such that trajectories starting in $W_1^{k_0}$ are mapped to $\kappa_{12}^{k_0}\left( \Sigma_1^{\textnormal{out}, k_0}\right) =  \Sigma_2^{\textnormal{in}, k_0}$.
\begin{proposition} \label{propK2}
Fix $k_0 \geq 2$. Then the following results hold for small $\delta >0$ depending on $\mu$:
\begin{enumerate}
\item[(P1)] If $\mu < 1$, every trajectory starting in $\Sigma_2^{\textnormal{in}, k_0}$ 
passes through $\Sigma_{2,\txta }^{\textnormal{out}, k_0}$, and, hence, so does $M_{\txta,2}^{-, k_0}$.
\item[(P2)] If $\mu > 1$, every trajectory starting in $\Sigma_2^{\textnormal{in}, k_0}$ 
passes through  $\Sigma_{2,\txte }^{\textnormal{out}, k_0}$ and, hence, so does $M_{\txta,2}^{-, k_0}$.
\end{enumerate}
\end{proposition}
\begin{proof}
%
Set $ K := \max \{C_{u, k_0}, C_{v, k_0}\} \sum_{k=2}^{k_0} \frac{1}{\left|b_k\right|}$. Then we have, at $\Sigma_2^{\textnormal{in}, k_0}$
\be
\label{eq:sum_estimate}
\frac{1}{2 a_2} \sum_{j=2}^{k_0} \left(u_{j,2}^2 - v_{j,2}^2 \right) \leq \frac{\delta^{-1/4}}{2 \nu} K \delta^{3/4} < \left| \mu -1 \right|,
\ee
when
$$ \delta < \frac{4 \nu^2}{K^2}  \left| \mu -1 \right|^2.$$
Hence, the sign of $(\mu-1)$ does not change and the dynamics for $u_{1,2}, v_{1,2}$ are as in the standard transcritical ODE case for sufficiently small $r_1$, as long as $v_{k,2}, u_{k,2}$, $k \geq 2$, do not grow thereafter.

The global stability of $v_{k,2}=0$ for all $k \geq 2$ is clear from equation~\eqref{eq:k0_K_2_general_vk} for any $r_2 > 0$ (and $a_2$ small enough in comparison to potential higher order terms $\mathcal{O} \left( r_2^8\right)$), or  $v_{k,2}$ staying constant for $r_2=0$ respectively.
For sufficiently small $r_2$, the linear stability of $u_{k,2} =0$ --- and due to~\eqref{eq:sum_estimate} this is sufficient for the neighbourhoods we are considering --- only becomes an issue for large $u_{1,2}$, as can be seen from~\eqref{eq:u12_bound}; but even in this case, by applying the map $\kappa_{21}^{k_0}$~\eqref{kappa21} 
and taking into account that $u_{1,2} \leq \delta^{-3/8}$, we observe that equation~\eqref{eq:u12_bound} is equivalent to
$$\nu  < \frac{\pi^2}{\sqrt{2}},$$
which is, of course, easily satisfied by choosing $\nu$ accordingly.

These considerations show that we can apply the same reasoning as in the proof of~\cite[Proposition 3.6]{EngelKuehn} for $h \to 0$, or the corresponding proof in \cite{KruSzm4}, to obtain the corresponding behaviour for $u_{1,2}, v_{1,2}$ and by that the claim. For the precise scaling of $v_{1,2}$ in $\Sigma_{2,\txte }^{\textnormal{out}, k_0}$, we refer to the detailed proof of~\cite[Proposition 3.6]{EngelKuehn}.
\end{proof}
As for chart $K_1$, we can give the corresponding PDE to the Galerkin ODEs. In chart $K_2$, this obviously only concerns a rescaling such that the PDE associated to system~\eqref{eq:k0_K_2_general} reads simply, for any $r_2 \geq 0$,
\be
\label{eq:mainPDE_for_calculations_K2}
\begin{array}{lcl}
 \partial_{\tau} u_2&=&\sqrt{2a_2} \left(\partial_x^2 u_2 + u^2-v^2 + \mu + \mathcal{O}\left(r_2^3\right)\right),\\
\partial_{\tau} v_2 &=& \sqrt{2a_2}\left( r_2^8 \partial_x^2 v + 1 + \mathcal{O}\left(r_2^{8} \right) \right), \ u_2=u_2(x,t),~v_2=v_2(x,t),~(x,t)\in[-a_2,a_2]\times[0,T],
\end{array}
\ee
under changing time back via $t'=\sqrt{2a_2}t$. 
%
\subsection{Analysis in third chart}
Finally, we consider the additional \emph{exit chart} $K_3$, determined by setting
$\bar u_1 = 1$. The coordinate transformation is given by
\begin{equation} \label{K3_scaling}
u_1 = r_3^{3} , \ v_1 = r_3^{3}v_{1,3},  \ u_k = r_3^{3} u_{k,3}, \ v_k = r_3^{3} v_{k,3}, \, k\geq 2, \
\varepsilon = r_3^8 \varepsilon_3, \ a = r_3^{-2} a_3,
\end{equation}
where $r_3 \geq 0$.

We can derive the following equations for the dynamics in chart $K_3$:
\begin{proposition} \label{prop:K3}
Consider the rescaling~\eqref{K1_scaling}.
Under this change of coordinates, system~\eqref{eq:k0_general}, for any $k_0 \in \mathbb{N}$, becomes 
\begin{subequations}
\label{eq:k0_K_3_general}
\begin{align}
\partial_{\tau} r_3&= \frac{1}{3} r_3 F_3(r_3, v_{1,3}, \varepsilon_3, a_3, u_{2,3}, v_{2,3}, \dots)\,, \label{eq:r_3}\\
\partial_{\tau} \varepsilon_3&= -\frac{8}{3} \varepsilon_3 F_3(r_3, v_{1,3}, \varepsilon_3, a_3, u_{2,3}, v_{2,3}, \dots)\,, \\
\partial_{\tau} v_{1,3}&= 2a_3\varepsilon_3  - v_{1,3} F_3(r_3, v_{1,3}, \varepsilon_3, a_3, u_{2,3}, v_{2,3}, \dots) + \mathcal O \left(r_3^{8}\right),\\
\partial_{\tau} u_{k,3}&= b_k a_3^{-3/2} u_{k,3} + 2 (u_{k,3} - v_{k,3}v_{1,3}) 
+ \sqrt{2}\sum_{i,j=2}^{k_0} \alpha_{i,j}^k (u_{i,3} u_{j,3} - v_{i,3} v_{j,3}) \nonumber\\
&-  u_{k,3}F_3(r_3, v_{1,3}, \varepsilon_3, a_3, u_{2,3}, v_{2,3}, \dots) + \mathcal O \left(r_3^3\right), \ k \geq 2, \\
\partial_{\tau} v_{k,3} &=  a_3^{-3/2} \left( r_3^{8}  b_k  v_{k,3} \varepsilon_3  -  a_3^{3/2} v_{k,3}F_3(r_3, v_{1,3}, \varepsilon_3, a_3, u_{2,3}, v_{2,3}, \dots) + a_3^{3/2}\mathcal O \left(r_3^{8}\right) \right), \ k \geq 2,\\
\partial_{\tau} a_3&= \frac{2}{3} a_3 F_3(r_3, v_{1,3}, \varepsilon_3, a_3, u_{2,3}, v_{2,3}, \dots)\,,
\end{align} 
\end{subequations}
where 
\begin{equation}
F_3(r_3, v_{1,3}, \varepsilon_3, a_3, u_{2,3}, v_{2,3}, \dots) = (1 - v_{1,3}^2) + 2a_3 \varepsilon_3 \mu + \sum_{j=2}^{k_0} u_{j,3}^2 - v_{j,3}^2  + \mathcal O \left(r_3^3\right).
\end{equation}
\end{proposition}
\begin{proof}
Similarly to proof of Proposition~\ref{prop:K1}.
\end{proof}
The dynamics are organized similarly to the dynamics in the chart $K_1$ with permuted stability. In more detail, for any $k_0 \geq 1$, note the existence of the three groups of equilibria
\begin{align} 
P_{\textnormal{r},3}^{-, k_0} &= \left\{ a_3 > 0\,:\,p_{\textnormal{r},3}^{-, k_0}(a_3) = \left(0,0,-1,0,0, \dots, a_3\right) \right\}, \label{equilibria_K3_ar_min} \\
P_{\textnormal{r},3}^{+, k_0} &=\left\{ a_3 > 0\,:\,p_{\textnormal{a},3}^{+, k_0}(a_3) = \left(0,0,1,0,0, \dots, a_3\right) \right\},\label{equilibria_K3_ar_plus}
\end{align}
and 
\begin{equation} \label{equilibria_K3_neutral}
Q_{3}^{\textnormal{out}, k_0} = \left\{a_3 > 0\,: \, q_{3}^{\textnormal{out}, k_0}(a_3) = (0,0,0,0,0 \dots,  a_3)\right\}.
\end{equation}
The equilibrium $q_{3}^{\textnormal{out}, k_0}(a_3)$ is a saddle with eigenvalue $-1$ in the $v_{1,1}$-direction.

\begin{lemma} \label{lem:K3_equilibira}
The equilibria $p_{\textnormal{r},3}^{-, k_0}(a_3)$ and $p_{\textnormal{r},3}^{+, k_0}(a_3)$ are 
\begin{enumerate}
\item unstable in the $v_{1,3}$-direction with eigenvalue $2$,
\item stable in the $u_{k,3}$-directions with eigenvalues $2+ b_k a_3^{-3/2}$ (for $a_3$ not too large),
\item  and neutral in all other directions.
\end{enumerate}
Furthermore, for $0 < r_3 < R_3$ sufficiently small, there are normally hyperbolic sets of equilibria $S_{\txtr, 3}^{-, k_0}$ and $S_{\txtr, 3}^{+, k_0}$ emanating from the lines $P_{\textnormal{r},3}^{-, k_0}$ and $P_{\textnormal{r},3}^{+, k_0}$.
\end{lemma}
\begin{proof}
This is analogous to Lemma~\ref{lem:K1_equilibira}.
\end{proof}
We now need to investigate the transition from $K_2$ to $K_3$. To this end, we compute:
\begin{lemma} \label{lem:K3K2}
The changes of coordinates between the charts $K_2$ and $K_3$ for $k_0 \in \mathbb{N}$ are given by the maps
$\kappa_{32}^{k_0}: K_3 \to K_2$ with
\begin{equation} \label{kappa32}
u_{1,2} = \varepsilon_{3}^{-3/8}, \ v_{1,2} = \varepsilon_{3}^{-3/8} v_{1,3}, \ u_{k,2} = \varepsilon_3^{-3/8} u_{k,3}, \ v_{k,2} = \varepsilon_3^{-3/8} v_{k,3}, \ r_2 = \varepsilon_3^{1/8} r_3, \ a_2 = \varepsilon_3^{1/4} a_3,
\end{equation}
and
$\kappa_{23}^{k_0}: K_2 \to K_3$ with
\begin{equation} \label{kappa23}
\varepsilon_{3} = u_{1,2}^{-8/3}, \ v_{1,3} = u_{1,2}^{-1} v_{1,2}, \ u_{k,3} = u_{1,2}^{-1} u_{k,2}, \ v_{k,3} = u_{1,2}^{-1} v_{k,2}, \ r_3 = u_{1,2}^{1/3} r_2, \ a_3 = u_{1,2}^{2/3} a_2.
\end{equation}
\end{lemma}
\begin{proof}
Follows from straight-forward calculations.
\end{proof}
In accordance with the discussion before Proposition~\ref{propK2}, we set
$$\Sigma_3^{\textnormal{in}, k_0} = \left\{ \varepsilon_3 = \delta \right\},$$
and 
$$ \Sigma_3^{\textnormal{out}, k_0} = \{ r_3 = \rho \}.$$
\begin{lemma} \label{lem:trans_K3}
The transition map $\Pi_3$  from $\Sigma_3^{\textnormal{in}, k_0}$ to $ \Sigma_3^{\textnormal{out}, k_0}$ is well-defined on $\kappa_{23}^{k_0}\left(\Sigma_{2,\txte }^{\textnormal{out}, k_0} \right)$. Furthermore, for any $z \in \kappa_{23}^{k_0}\left(\Sigma_{2,\txte }^{\textnormal{out}, k_0} \right) \subset \Sigma_3^{\textnormal{in}, k_0}$, we have that $\pi_{v_{1,3}}\left(\Pi_3(z)\right) = \mathcal{O}\left( \delta^{1/4} \right)$.
\end{lemma}
\begin{proof}
Similarly to \cite[Proposition 3.9]{EngelKuehn}. The second statement follows directly from using $v_{1,3} = u_{1,2}^{-1} v_{1,2}$.
\end{proof} 
\section{Main result with proof}
\label{sec:mainresult}
\begin{figure}[ht]
        \centering
        \begin{subfigure}{0.5\textwidth}
        \centering
  		\begin{overpic}[width=1\linewidth]{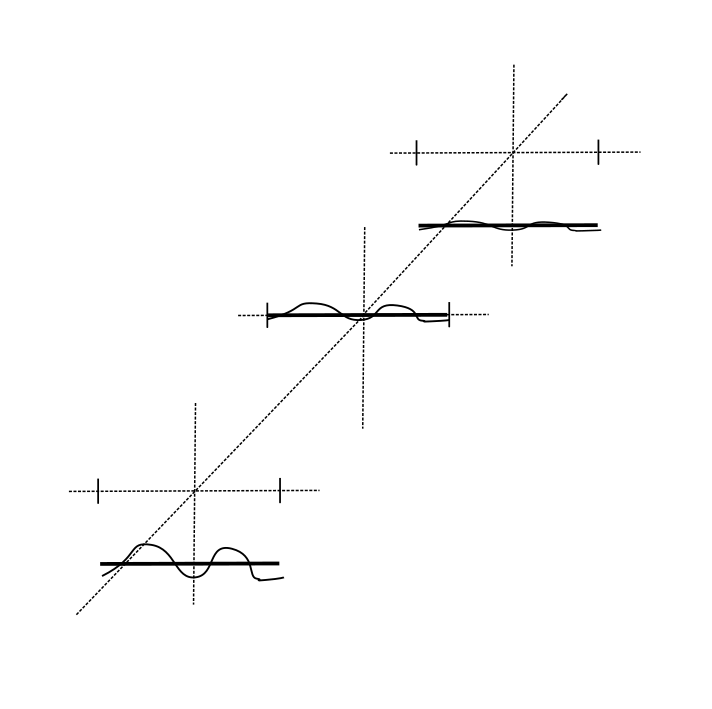}   
\put(45,20){\scriptsize $u^*=-\rho$}
  	\put(70,55){\scriptsize $u^*=0$}	
  		\put(86,68){\scriptsize $u^*=-\rho$}	
	\put(21,40){\tiny $u(x)$}
  		\put(10,27){\tiny $-a$}
  		\put(38,27){\tiny $a$}
  		\put(30,77){\tiny $\sgn(v_1) \|v\|_2 \approx \rho$}
  	        
                \put(10,55){\tiny $\sgn(v_1) \|v\|_2 =0$}
                \put(50,30){\tiny $\sgn(v_1) \|v\|_2 \approx - \rho$}
        \end{overpic}
        \caption{$\mu <1$}
        \label{musmaller1}
		\end{subfigure}%
        \begin{subfigure}{0.5\textwidth}
        \centering
  		\begin{overpic}[width=1\linewidth]{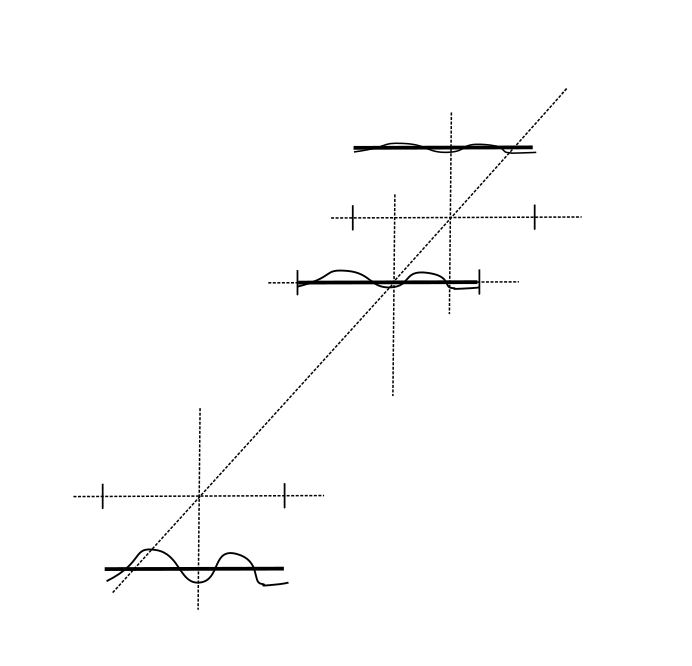}   
\put(45,13){\scriptsize $u^*=-\rho$}
  	\put(78,55){\scriptsize $u^*=0$}	
  		\put(85,74){\scriptsize $u^*=\rho$}	
	\put(21,35){\tiny $u(x)$}
  		\put(11,21){\tiny $-a$}
  		\put(40,21){\tiny $a$}
  		\put(15,65){\tiny $\sgn(v_1) \|v\|_2 \approx \mathcal O (\epsilon^{1/4} )$}
  	        
                \put(11,54){\tiny $\sgn(v_1) \|v\|_2 =0$}
                \put(52,24){\tiny $\sgn(v_1) \|v\|_2 \approx - \rho$}
        \end{overpic}
        \caption{$\mu >1$}
        \label{mularger1}
		\end{subfigure}
		\caption{Sketch of the dynamics of a typical solution of the $k_0$-Galerkin truncation of the PDE system~\eqref{eq:mainPDE_fasttime} along a slow attracting manifold close to the critical manifold as described in Theorem~A for $\mu < 1$ (a) and $\mu > 1$ (b).}
        \label{fig:main_theorem}
\end{figure}

For sufficiently small $\rho > 0$ and $k_0 \in \mathbb N$, we introduce
\begin{itemize}
\item $ \Delta^{\textnormal{in}, k_0}$ as a small neighbourhood around the pair of constant functions $$(u^*, v^*) \equiv (-\rho, -\rho) \in L^2([-a,a]) \times L^2([-a,a]),$$
where the Fourier coefficients for $k > k_0$ are all zero,
\item $\Delta_\txta^{\textnormal{out}, k_0}$ as a small neighbourhood around the pair of constant functions $$(u^*, v^*) \equiv (-\rho, \rho) \in L^2([-a,a]) \times L^2([-a,a]),$$
where the Fourier coefficients for $k > k_0$ are all zero,
\item $\Delta_\txte^{\textnormal{out}, k_0}$ as a small neighbourhood around the pair of constant functions $$(u^*, v^*) \equiv (\rho, 0) \in L^2([-a,a]) \times L^2([-a,a]),$$
where the Fourier coefficients for $k > k_0$ are all zero.
\end{itemize}
Furthermore, let us denote by $\Pi_\txta$ and $\Pi_\txte$ the transition maps along the dynamics of  system~\eqref{eq:mainPDE_fasttime} from $ \Delta^{\textnormal{in}, k_0}$ to $\Delta_\txta^{\textnormal{out}, k_0}$ and $\Delta_\txte^{\textnormal{out}, k_0}$ respectively. We can now give the following detailed version of our main result which is illustrated in Figure~\ref{fig:main_theorem}.
\begin{theorem}[Detailed Version of Theorem~\ref{thm:mainresult}] 
\label{thm:mainresult_details}
The attracting slow manifolds $S_{\txta, \varepsilon}^{-,k_0}$ and $S_{\txta, \varepsilon}^{+,k_0}$ near the origin for system~\eqref{eq:mainPDE_Galerkin}, truncated at $k_0 \in \mathbb{N}$,
exhibit the following behavior:

For any fixed, $k_0 \in \mathbb N$, $\mu \neq 1$ and sufficiently small $\rho > 0$, there is an $\varepsilon = \varepsilon_0(k_0) > 0$ such that for all $\varepsilon \in (0, \varepsilon_0]$:
\begin{enumerate}
\item[(T1)] If $ \mu < 1$, then the neighbourhood $\Delta^{\textnormal{in}, k_0}$ (including $\Delta^{\textnormal{in}, k_0} \cap S_{\txta, \varepsilon}^{-,k_0}$) is mapped by $\Pi_\txta$ 
to a set containing $S_{\txta, \varepsilon}^{+, k_0}$.
\item[(T2)] If $ \mu > 1$, then the manifold $S_{\txta, \varepsilon}^{-, k_0}$ passes through 
$\Delta_\txte^{\textnormal{out}, k_0}$ in close vicinity to a point $(\rho, p(\varepsilon))$ where $p(\varepsilon) = 
\mathcal{O}(\varepsilon^{1/4})$. The section $\Delta^{\textnormal{in}, k_0}$ is mapped by $\Pi_\txte$ 
to a set containing $S_{\txta, \varepsilon}^{-, k_0} \cap \Delta_\txte^{\textnormal{out}, k_0}$.
\end{enumerate}
\end{theorem}
\begin{proof}
(T1) and (T2) can be now proven putting everything in Section~\ref{sec:blowup} together.
In more details, we firstly translate, for $\varepsilon_0$ and $\rho$ small enough, the set $ \Delta^{\textnormal{in}, k_0}$ into a set $W_1^{k_0} \subset \Sigma_1^{\textnormal{in}, k_0}$ which forms a neighbourhood of $M_{\txta,1}^{-, k_0} \cap \Sigma_1^{\textnormal{in}, k_0}$ for sufficiently small $v_{k,1}$ such that trajectories starting in $W_1^{k_0}$ are mapped to $\kappa_{12}^{k_0}\left( \Sigma_1^{\textnormal{out}, k_0}\right) =  \Sigma_2^{\textnormal{in}, k_0}$.
According to Proposition~\ref{propK2}, we have that, if $\mu < 1$, every trajectory starting in $\Sigma_2^{\textnormal{in}, k_0}$ 
passes through $\Sigma_{2,\txta }^{\textnormal{out}, k_0}$, and so does $M_{\txta,2}^{-, k_0}$. It is now easy to see, with a reversed calculation to the transition from $K_1$ to $K_2$, that $ \kappa_{21}^{k_0}\left(\Sigma_{2,\txta }^{\textnormal{out}, k_0}\right)$, with $a_1 = \nu$, $\varepsilon_1 = \delta$, intersects with $ M_{\txta,1}^{-, k_0}$ and $M_{\txta,1}^{+, k_0}$. Claim (T1) follows then by tracking the dynamics in chart $K_1$ up to reaching, in original coordinates, $\Delta_\txta^{\textnormal{out}, k_0}$.

If $\mu > 1$, every trajectory starting in $\Sigma_2^{\textnormal{in}, k_0}$ 
passes through  $\Sigma_{2,\txte }^{\textnormal{out}, k_0}$ and so does $M_{\txta,2}^{-, k_0}$. Analyzing the dynamics passing through $\kappa_{23}^{k_0}\left(\Sigma_{2,\txte }^{\textnormal{out}, k_0} \right)$ in chart $K_3$, Lemma~\ref{lem:trans_K3} gives the result (T2) under transformation back to original coordinates.
\end{proof}
We remark that the case $\mu =1$ coincides with the phenomenon of \textit{canards} which requires a whole new study for the situation of PDEs. 

\section{Conclusion \& Outlook}
\label{sec:conclusion}

In this work, we have provided a blow-up approach to a fast-slow partial differential equation in the context of a dynamic transcritical bifurcation for reaction-diffusion systems. We have employed a Galerkin discretization, made the domain length a dynamic variable, and then performed a blow-up on each finite-dimensional truncation level. In three different charts, we have controlled the invariant manifolds, which allowed us to track the attracting slow manifold through the fast-slow transcritical point, where normal hyperbolicity is lost. Furthermore, we identified desingularized PDEs in all three charts. Effectively, our main result is a fully geometric version of a dynamic Crandall-Rabinowitz theorem, in arbitrarily close finite-dimensional approximation.

    Note that we have focused on Neumann boundary conditions for two main reasons. Firstly, this is a natural boundary assumption for many application-oriented models in fluid dynamics, thermodynamics or also spatial ecological models. Secondly, from a mathematical point of view, these boundary conditions allow for anchoring the behaviour around non-trivial homogeneous solutions that enable an analysis that is easily comparable with low-dimensional ODE models. In particular, for having a first go on the complicated geometric blow-up problem of PDEs and its Galerkin approximations, this implication of Neumann boundary conditions has proven very helpful, as seen in this paper. Dirichlet boundary conditions, for example, already entail a more delicate situation, without non-trivial homogeneous solutions and possible effects like boundary layers that have to be included in the analysis of passing through a singularity.

Having a blueprint of blowing up infinite-dimensional multiscale dynamical systems, one can now aim to tackle many other related issues. 
One may aim to generalize our current transcritical setting in different ways, e.g., to very small or very large domains, to Dirichlet/mixed or even singular boundary data, to singular diffusion constants, to quasi-linear parabolic equations, to canard solutions, and to very general single eigenvalue crossing problems involving a trivial homogeneous branch such as pitchfork singularities. Other low-codimension dynamic/fast-slow bifurcation points for reaction diffusion PDEs would also be of interest including higher-order folded singularities (as an extension of~\cite{Engeletal2022} or Hopf points). In particular, these results should immediately be helpful to understand pattern formation dynamics, amplitude/modulation equations with slowly-varying parameters, and to make many results obtained by formal asymptotic matching for PDEs rigorous via a geometric approach. From the viewpoint of PDE theory, it seems also natural to conjecture that upon a suitable spatial discretization, one may be able to treat other classes of differential operators such as wave-type operators, advective terms, Fokker-Planck-type equations, and similar classes of problems, where other dynamical systems methods have proven to be successful already. From a methodological viewpoint, exploring the connection to free boundary problems and PDEs on Banach manifolds, is a natural continuation of our work.\medskip 

\textbf{Acknowledgments:} The authors gratefully acknowledge support by the DFG via the SFB TRR 109 Discretization in Geometry and Dynamics. ME additionally thanks the DFG for support via Germany's Excellence Strategy -- The Berlin Mathematics Research Center MATH+ (EXC-2046/1, project ID: 390685689). CK likes to also thank the VolkswagenStiftung for partial support via a Lichtenberg Professorship. CK also acknowledges interesting discussions with Peter Szmolyan in 2010/2011 (when CK started working on the problem of extending the blow-up method to PDEs), as well as more recent general discussions regarding Allen-Cahn-type PDEs with Daniel Matthes. CK and ME want to thank Felix Hummel for several insightful comments and discussions regarding an early version of this work.

\appendix

\section{Appendix}

\subsection{Truncation at third mode}
\label{ap:truncate3}
To provide an even better grasp of the structure of the Galerkin equations, we also considered $k_0 =3$ for the computation of the center/slow manifold before the bifurcation point: in this case equation~\eqref{eq:k0_general} reads
\begin{equation} \label{eq:k03}
\begin{array}{lcl}
\partial_t u_1&=& u_1^2 - v_1^2 + 2a \varepsilon \mu +  u_2^2 - v_2^2 + u_3^2 - v_3^2,\\
\partial_t v_1&=& 2a\varepsilon , \\
\partial_t u_2&=& \hat \lambda_2 u_2 + 2 (u_2 u_1  -  v_2 v_1)+  \sqrt{2}(u_2 u_3  -  v_2 v_3),\\
\partial_t v_2&=& \varepsilon \hat \lambda_2 v_2, \\
\partial_t u_3&=& \hat \lambda_3 u_3 + 2 (u_3 u_1  -  v_3 v_1)+  \frac{1}{\sqrt{2}}(u_2^2 -  v_2^2),\\
\partial_t v_3&=& \varepsilon \hat \lambda_3 v_3, \\
\partial_t \varepsilon &=& 0,
\end{array}
\end{equation} 
or again, with shift to $\tilde u_1, \tilde v_1$,
\begin{equation} \label{eq:k03_shifted}
\begin{array}{lcl}
\partial_t \tilde u_1&=& (\tilde u_1+c)^2 - (\tilde v_1+c)^2 + 2a \varepsilon \mu +  u_2^2 - v_2^2 ,\\
\partial_t \tilde v_1&=& 2a\varepsilon , \\
\partial_t u_2&=& \hat \lambda_2 u_2 + 2 (u_2 (\tilde u_1 + c)  -  v_2 (\tilde v_1 + c)) +  \sqrt{2}(u_2 u_3  -  v_2 v_3),\\
\partial_t v_2&=& \varepsilon \hat \lambda_2 v_2, \\
\partial_t u_3&=& \hat \lambda_3 u_3 + 2 (u_3 (\tilde u_1 + c)  -  v_3 (\tilde v_1 + c))+  \frac{1}{\sqrt{2}}(u_2^2 -  v_2^2),\\
\partial_t v_3&=& \varepsilon \hat \lambda_3 v_3, \\
\partial_t \varepsilon &=& 0.
\end{array}
\end{equation}
Again considering the equilibrium at the origin, we  obtain, compared to the case $k_0 =2$, the additional stable direction in $u_3$ with eigenvalue $\hat \lambda_3 + 2c < 0$ 
and center direction in $v_2$. Proceeding analogously to the case $k_0=2$ and transferring to standard form coordinates
\begin{align*}
\tilde u_1 &= y_1 + v_1 - 2 a \varepsilon \left( \frac{\mu -1}{2c} \right), \\
u_2 &= y_2 + \frac{2c}{2c + \hat \lambda_2}v_2,\\
u_3 &= y_3 + \frac{2c}{2c + \hat \lambda_3}v_3,\\
\tilde v_1 &= x_1,\
\varepsilon = \frac{x_2}{2a},\
v_2 = x_3,\
v_3 = x_4,
\end{align*}
we consider
\begin{align*}
\hat h_1(x_1, x_2, x_3, x_4) &= b_{11} x_1^2 + b_{12} x_1 x_2 + b_{22} x_2^2 + b_{13} x_1 x_3 + b_{23} x_2 x_3 + b_{33} x_3^2 \\
&+b_{44} x_4^2 + b_{14} x_1 x_4 + b_{24} x_2 x_4 + b_{34} x_4 x_3  + \mathcal{O}(3), \\
\hat h_2(x_1, x_2, x_3, x_4) &= c_{11} x_1^2 + c_{12} x_1 x_2 + c_{22} x_2^2 + c_{13} x_1 x_3 + c_{23} x_2 x_3 + c_{33} x_3^2 \\
&+c_{44} x_4^2 + c_{14} x_1 x_4 + c_{24} x_2 x_4 + c_{34} x_4 x_3  + \mathcal{O}(3), \\
\hat h_3(x_1, x_2, x_3, x_4) &= d_{11} x_1^2 + d_{12} x_1 x_2 + d_{22} x_2^2 + d_{13} x_1 x_3 + d_{23} x_2 x_3 + d_{33} x_3^2 \\
&+d_{44} x_4^2 + d_{14} x_1 x_4 + d_{24} x_2 x_4 + d_{34} x_4 x_3  + \mathcal{O}(3).
\end{align*}
The associated invariance equation gives the coefficients,
\begin{align} 
b_{11} &= b_{13}= b_{23}=b_{14}=b{24}= b_{34} =0, \ b_{12} = \frac{\mu -1}{2c^2}, \ b_{22} = \frac{(\mu-3)(\mu -1)}{8c^3}, \nonumber \\ b_{33} &= \frac{1}{2c} - \frac{2c}{(\hat \lambda_2 + 2c)^2}, \ b_{44}= \frac{1}{2c} - \frac{2c}{(\hat \lambda_3 + 2c)^2}, \label{eq:coeff_k03_b}\\
c_{11} &=  c_{22} = c_{33}= c_{44} = c_{12}= c_{14}=c_{24}= 0, \ c_{13} = \frac{2 \hat \lambda_2}{(\hat \lambda_2 + 2c)^2}, \nonumber\\
c_{23} &= \frac{c \hat \lambda_2(2c+\hat \lambda_2)+4 ac(\mu-1)+2a \hat \lambda_2 \mu}{a(\hat \lambda_2 + 2c)^3}, \ c_{34}= \frac{\sqrt{2}( \hat \lambda_2 \hat \lambda_3 + 2c(\hat \lambda_2 + \hat \lambda_3))}{(2c+\hat \lambda_2)^2(2c+\hat \lambda_3)},\\
d_{11}&= d_{22}=d_{44} =d_{12}=d_{13}=d_{34}=d_{23}= =0, \ d_{13} = \frac{2 \hat \lambda_3}{(\hat \lambda_3 + 2c)^2}, \nonumber \\
d_{24} &= \frac{c \hat \lambda_3(2c+\hat \lambda_3)+4 ac(\mu-1)+2a \hat \lambda_3 \mu}{a(\hat \lambda_3 + 2c)^3}, \ d_{33}= \frac{ \hat \lambda_2 (4c + \hat \lambda_2)}{\sqrt{2}(2c+\hat \lambda_2)^2(2c+\hat \lambda_3)}.
\end{align}
Transforming back gives the center manifold approximation in form of the maps
\begin{align}
\tilde u_1 &= \tilde h_1 (\tilde v_1, \varepsilon, v_2, v_3) =\tilde v_1 - \frac{a(\mu -1)}{c} \varepsilon + \frac{a(\mu-1)}{c^2}\tilde v_1 \varepsilon \nonumber\\
&- \frac{a^2(\mu-3)(\mu -1)}{2c^3} \varepsilon^2 + \left( \frac{1}{2c} - \frac{2c}{(\hat \lambda_2 + 2c)^2}\right) v_2^2 + \left(\frac{1}{2c} - \frac{2c}{(\hat \lambda_3 + 2c)^2}\right) v_3^2 + \mathcal{O}(3), \label{eq:cmtilde_k03_1} \\
\tilde u_2 &= \tilde h_2 (\tilde v_1, \varepsilon, v_2, v_3) = \frac{2c}{2c + \hat \lambda_2} v_2 +  2 \frac{c \hat \lambda_2(2c+\hat \lambda_2)+4 ac(\mu-1)+2a \hat \lambda_2 \mu}{(\hat \lambda_2 + 2c)^3} v_2 \varepsilon \nonumber \\
&+ \frac{2 \hat \lambda_2}{(\hat \lambda_2 + 2c)^2} \tilde v_1 v_2 + \frac{\sqrt{2}( \hat \lambda_2 \hat \lambda_3 + 2c(\hat \lambda_2 + \hat \lambda_3))}{(2c+\hat \lambda_2)^2(2c+\hat \lambda_3)} v_2 v_3 +  \mathcal{O}(3), \label{eq:cmtilde_k03_2} \\
\tilde u_3 &= \tilde h_3 (\tilde v_1, \varepsilon, v_2, v_3) = \frac{2c}{2c + \hat \lambda_3} v_3 +  2 \frac{c \hat \lambda_2(2c+\hat \lambda_3)+4 ac(\mu-1)+2a \hat \lambda_3 \mu}{(\hat \lambda_3 + 2c)^3} v_3 \varepsilon \nonumber \\
&+ \frac{2 \hat \lambda_3}{(\hat \lambda_3 + 2c)^2} \tilde v_1 v_3 + \frac{ \hat \lambda_2 (4c + \hat \lambda_2)}{\sqrt{2}(2c+\hat \lambda_2)^2(2c+\hat \lambda_3)} v_2^2 +  \mathcal{O}(3). \label{eq:cmtilde_k03_3}
\end{align}
One can check by computation that, indeed, $\tilde h= (\tilde h_1, \tilde h_2, \tilde h_3)$ satisfies the system of invariance equations associated with~\eqref{eq:k03_shifted}, given by
\begin{equation} \label{eq:invar_k03_shifted}
\begin{array}{lcl}
 2 a \varepsilon \partial_{\tilde v_1} \tilde h_1 + \varepsilon \hat \lambda_2 v_2 \partial_{\tilde v_2} \tilde h_1 + \varepsilon \hat \lambda_3 v_3 \partial_{\tilde v_3} \tilde h_1  &=& (\tilde h_1 + c)^2 - (\tilde v_1 +c)^2 + 2 a \varepsilon \mu \\
 &+& (\tilde h_2)^2 - v_2^2 + (\tilde h_3)^2 - v_3^2 + \mathcal{O}(3),\\
 2 a \varepsilon \partial_{\tilde v_1} \tilde h_2 + \varepsilon \hat \lambda_2 v_2 \partial_{\tilde v_2} \tilde h_2 + \varepsilon \hat \lambda_3 v_3 \partial_{\tilde v_3} \tilde h_2 &=& \hat \lambda_2 \tilde h_2 + 2 (\tilde h_2 (\tilde h_1 + c)  -  v_2 (\tilde v_1 + c))\\
 &+&  \sqrt{2}(\tilde h_2 \tilde h_3  -  v_2 v_3) + \mathcal{O}(3),\\
  2 a \varepsilon \partial_{\tilde v_1} \tilde h_3 + \varepsilon \hat \lambda_2 v_2 \partial_{\tilde v_2} \tilde h_3 + \varepsilon \hat \lambda_3 v_3 \partial_{\tilde v_3} \tilde h_3 &=& \hat \lambda_3 \tilde h_3 + 2 (\tilde h_3 (\tilde h_3 + c)  -  v_3 (\tilde v_1 + c))\\
 &+&  \frac{1}{\sqrt{2}}\left(\tilde h_2^2  -  v_2^2 \right) + \mathcal{O}(3).
 \end{array}
\end{equation}
The relations $u_1 = h_1 ( v_1, \varepsilon, v_2, v_3)$, $u_2 = h_2 ( v_1, \varepsilon, v_2, v_3)$ and $u_3 = h_3 ( v_1, \varepsilon, v_2, v_3)$ are obtained by transforming back into original coordinates.
\subsection{Sketch: Alternative PDE argument for second chart}
\label{ap:alternative_proof}

We decided to work with Galerkin systems throughout our analysis to keep one consistent framework. Instead of working directly with the Galerkin systems, we want to mention that it is also possible to use more PDE-based arguments \emph{in each chart}. In particular, the second chart serves as a good illustration of this strategy, which we shall sketch here. Hence, we again go back to the PDE view for understanding the dynamics through chart $K_2$:
\beann
\partial_t u &=& \partial_x^2 u + u^2 - v^2 +\mu\varepsilon + \cO(\cdots),\\
\partial_t v &=& \varepsilon \partial_x^2 v + \varepsilon + \cO(\cdots),
\eeann
where $x\in(-a,a)$ and homogeneous Neumann boundary conditions hold. A direct scaling such as 
\benn
u=\varepsilon^{1/2}U,\quad v=\varepsilon^{1/2}V,\quad x=\varepsilon^{-1/4}X,\quad t=\varepsilon^{-1/2}T
\eenn
gives the PDE
\begin{equation}
\label{eq:PDE_appendix}
\begin{array}{lcl}
\partial_T U &=& \partial_X^2 U + U^2 - V^2 +\mu + \varepsilon^{p_u}\cO(\cdots),\\
\partial_T V &=& \varepsilon \partial_X^2 V + 1 +\varepsilon^{p_v} \cO(\cdots),
\end{array}
\end{equation}
posed for $X\in (-a\varepsilon^{1/4},a\varepsilon^{1/4})$. This form is, indeed, equivalent to equation~\eqref{eq:mainPDE_for_calculations_K2}. Denoting solutions of~\eqref{eq:PDE_appendix} by $(U_\varepsilon(T,X),V_\varepsilon(T,X))$, with initial data that are sufficiently close to a constant, as given by the transfer from the entry chart $K_1$, we observe that in $L^2[-a,a]$, as $\varepsilon\ra 0$,
\be
\label{eq:convchart2}
(U_\varepsilon(T,X),V_\varepsilon(T,X))\ra (U_0(T),V_0(T)) 
\ee 
where $(U_0(T),V_0(T))$ solves the ODE 
\be
\label{eq:ODE_appendix}
\begin{array}{lcl}
\partial_T U_0 &=& U_0^2 - V_0^2 +\mu,\\
\partial_T V_0 &=& 1.
\end{array}
\ee
The convergence follows from the fact that the Galerkin coefficients are uniformly bounded --- in particular, avoiding any divergences in $\varepsilon$ ---, as shown around Proposition~\ref{propK2}, such that the $L^2$-norm does not blow up as $\varepsilon\ra 0$.
In other words, the dynamics on the sphere, as expressed in chart $K_2$, coincide exactly as $\varepsilon\ra 0$, i.e., as $r_2\ra 0$, with the ODE dynamics for~\eqref{eq:ODE_appendix}. And by continuity in $\varepsilon$ in $L^2$, it then suffices over a finite time scale to use the ODE approximation for the PDE for sufficiently small $\varepsilon >0$; in fact, considering a finite time scale is sufficient, since the data in the entry and exit charts, that have to be matched, are located on a tube around suitable constant solutions. 
Hence, one could prove the statements of Proposition~\ref{propK2} by subtracting the solutions of equation~\eqref{eq:PDE_appendix} and~\eqref{eq:ODE_appendix}, assuming that initial data are close to a constant with some $\varepsilon$-dependent error, and use the results in \cite{KruSzm4}, or \cite{EngelKuehn} respectively.

In summary, the combination of a shrinking domain, the Neumann conditions, the parabolic regularity and the uniformly bounded $L^2$-norm in $\varepsilon$ effectively enforce that solutions $(U_\varepsilon(T,X),V_\varepsilon(T,X))$ stay near constants for longer and longer times, and, hence, in the limit $\varepsilon \to 0$ they must match the ODE solutions, which are completely independent of $X$. Note very carefully that the blow-up method has very nicely eliminated the higher-order reaction-terms in the limit, while the PDE approximation argument can take care of the Laplacians in the classical scaling chart $K_2$. This last step only works in the absolute singular limit $\varepsilon=0$ and the blow-up method guarantees that we can extend this result to a neighbourhood of the sphere as we do not require it over infinite time scales in blown-up space --- which we would do without the entry or exit charts.

\bibliographystyle{siamplain}
\bibliography{bibfile_PDE_blowup}

\end{document}